\documentclass[a4paper,12pt]{article}
\usepackage{amsmath,amssymb,amsfonts}
\usepackage{graphicx}
\usepackage{epsfig}
\usepackage{amsfonts}
\usepackage{amssymb, latexsym, amsmath, pb-diagram}
\usepackage{pstricks-add}
\usepackage{multicol}
\usepackage{bm}
\usepackage{mathrsfs}
\usepackage{xy}
\usepackage{epic,eepic}
\xyoption{all}

\numberwithin{equation}{section}

\newcommand{\zz}{\mathbb{Z}}
\def\XX{\mathscr X}

\def\RR{\mathscr R}

\usepackage{amssymb, latexsym, amsmath, pb-diagram}

\begin{document}
\title{The Cohomology Algebra of Polyhedral Product Spaces}
\author {
Qibing Zheng\\School of Mathematical Science and LPMC, Nankai University\\
Tianjin 300071, China\\
zhengqb@nankai.edu.cn\footnote{Project Supported by Natural Science
Foundation of China, grant No. 11071125\newline
\hspace*{5.5mm}Key words and phrases: polyhedral product,
diagonal tensor product\newline
\hspace*{5.5mm}Mathematics subject classification: 55N10}}\maketitle

\input amssym.def
\newsymbol\leqslant 1336
\newsymbol\geqslant 133E
\baselineskip=20pt
\def\w{\widetilde}
\begin{abstract} In this paper, we compute the integral singular cohomology ring of
homology split polyhedral product spaces and the singular cohomology algebra over a field of
polyhedral product spaces.
As an application, we give two polyhedral product spaces
${\cal Z}(K;X_1,A_1)$ and ${\cal Z}(K;X_2,A_2)$ such that
the cohomology homomorphisms $i_k^*\colon H^*(X_k)$ $\to H^*(A_k)$
induced by the inclusions are the same,
but the cohomology rings of the two polyhedral product spaces are not isomorphic.
\end{abstract}

\hspace*{50mm}${\displaystyle{\bf Table\,\,of\,\,Contents}}$

Section 1\, Introduction

Section 2\, Character Coproduct

Section 3\, Diagonal Tensor Product

Section 4\, Homology and Cohomology Group

Section 5\, Diagonal Tensor Product of Algebras and Coalgebras

Section 6\, Cohomology Algebra

Section 7\, Restriction Products
\vspace{3mm}

\newtheorem{Definition}{Definition}[section]
\newtheorem{Theorem}{Theorem}[section]
\newtheorem{Lemma}{Lemma}[section]
\newtheorem{Example}{Example}[section]
\font\hua=eusm10 scaled\magstephalf

\section{Introduction}

\hspace*{5.5mm}A polyhedral product is a relatively new construction with its origins
in Toric Topology and therefore closely related to toric objects coming from
algebraic and symplectic geometry. Since its formal appearance in work of
Buchstaber and Panov \cite{BP} and Grbi\'{c} and Theriault \cite{GR} in 2004,
the theory (especially the homotopical characteristics) of polyhedral products
has been developing rapidly. Due to its combinatorial nature coming from the underlying
simplicial complex and being a product space, the polyhedral product functor was quickly
recognized  as a complex construction but at the same time approachable.
As a result, polyhedral products are nowadays used not only in topology and geometry
but also group theory (abstract and geometric) as a collection of spaces on which to test
major conjectures and to form a rather delicate insight in building new theories.

To a simplicial complex $K$ on $[m]$ and
a family of CW-pairs $(\underline{X},\underline{A})=\{(X_k,A_k)\}_{k=1}^m$,
the polyhedral product functor assigns a topological space\vspace{-2mm}
$${\cal Z}(K;\underline{X},\underline{A})=\cup_{\sigma\in K}D(\sigma)\subset X_1{\times}X_2{\times}{\cdots}{\times}X_m,\vspace{-2mm}$$
where $D(\sigma)\cong\prod_{i\in\sigma}X_i\times\prod_{j\not\in\sigma} A_j$.
Depending on the choice of $(\underline{X},\underline{A})$,
the space ${\cal Z}(K;\underline{X},\underline{A})$ has other names.
If all the CW-pairs $(X_k,A_k)$ are the same $(X,A)$,
the polyhedral product space is denoted by ${\cal Z}(K;X,A)$.
In case when $(X,A)=(D^2,S^1)$, ${\cal Z}(K;D^2,S^1)$ is called a moment-angle complex.
In this form, a moment-angle complex was first introduced by Buchstaber and Panov \cite{BP} and
was widely studied by mathematicians in the area of toric topology and geometry.
For example, a moment-angle complex is an equivariant deformation retract of the complement
space of the associated to $K$ complex coordinate subspace arrangement
which is related to the classifying space of colored braid groups and
the configuration space for different classical mechanical systems
(see \cite{A},\cite{BP},\cite{B},\cite{GM},\cite{GR},\cite{GJ},\cite{GS}).
Buchstaber and Panov \cite{BP} proved\vspace{-2mm}
$$H^*({\cal Z}(K;D^2,S^1);\Bbb Z)={\rm Tor}_{\Bbb Z[x_1,\cdots,x_m]}^*(\Bbb Z(K),\Bbb Z),\vspace{-2mm}$$
where $\Bbb Z(K)$ is the Stanley-Reisner ring of $K$.
This result connects toric topology to homological algebra and combinatorial commutative algebra.
Baskakov \cite {BS} and independently Denham and Suciu \cite{DE} studied the Massey product on  the cohomology algebra of moment-angle complexes.
${\cal Z}(K;\!C\!P^{^{\infty}}\!\!,*)$ is the Davis-Januszkiewicz space which can be also defined as the Borel construction
of a moment-angle complex (see \cite{DJ}), whose cohomology ring is the Stanley-Reisner ring of $K$.
${\cal Z}(K;I,1)$ is a triangulation of the ${\rm cubic}$ complex (see \cite{BP}).
There arises a need to unify algebraic, geometric and combinatorial construction
coming from simplicial complexes and product of spaces.

The unstable homotopy types of polyhedral product spaces were studied by
Grbi\'{c} and Theriault \cite {GR},\cite{GJ},\cite{GS} and Beben and Grbi\'{c} \cite{GB} while the stable
homotopy types were studied by
Bahri, Bendersky, Cohen and Gitler \cite{B1},\cite{B2},\cite{B3}and many others.
The singular cohomology ring $H^*({\cal Z}(K;\underline{X},\underline{A}))$ is not known
even if all $H^*(A_k)$, $H^*(X_k)$ and $i_k^*\colon H^*(X_k)\to H^*(A_k)$ are known.
In this paper, we first compute the integral (co)homology group of
homology split polyhedral product spaces and the (co)homology group over a field
of polyhedral product spaces in Theorem 4.8.
A homology split polyhedral product space ${\cal Z}(K;\underline{X},\underline{A})$
satisfies that
all $H_*(X_k)$, $H_*(A_k)$ and the kernel of $i_k\colon H_*(A_k)\to H_*(X_k)$ induced by inclusion
are free groups.
For such a homology split polyhedral product space,\vspace{-1mm}
$$H^*({\cal Z}(K;\underline{X},\underline{A}))\cong\oplus_{(\sigma,\omega)\in\XX_m}\,\,\,
\w H^{*-1}(K_{\sigma,\omega}){\otimes}H^*_{\sigma\!,\,\omega}(\underline{X},\underline{A}),\vspace{-1mm}$$
where $\XX_m=\{(\sigma,\omega)\,|\,\sigma,\omega\subset[m],\,\,\sigma{\cap}\omega=\emptyset\}$,
$K_{\sigma,\omega}=({\rm link}_{_K}\sigma)|_\omega$,
$\w H^*(-)$ means reduced simplicial cohomology and\vspace{-2mm}
$$H^*_{\sigma\!,\,\omega}(\underline{X},\underline{A})
=H^1{\otimes}{\cdots}{\otimes}H^m,\quad
H^k=\left\{\begin{array}{cl}
{\rm ker}\,i_k^*&{\rm if}\, k\in\sigma,\vspace{1mm}\\
{\rm coker}\,i_k^*&{\rm if}\, k\in\omega,\vspace{1mm}\\
{\rm im}\,i_k^*&{\rm otherwise}.
\end{array}\right.$$

The (co)homology group is computed from the point of view of diagonal tensor product
defined in Section 3.
The diagonal tensor product of groups are naturally generalized to (co)algebras in Section 5.
In Theorem 6.9, we prove that the above diagonal tensor product of groups is a diagonal tensor product
of algebras\vspace{-2mm}
$$\big(H^*({\cal Z}(K;\underline{X},\underline{A})),\cup\big)
\cong \big(H_{\XX_m}^*(K)\otimes_{\XX_m} H_{\XX_m}^*(\underline{X},\underline{A}),\,
\cup_{K}\,{\otimes}_{\XX_m}\,\pi_{(\underline{X},\underline{A})}\big).\vspace{-2mm}$$
This result includes all the cases mentioned above and
their cohomology ring is computed in detail in Section 7.
The ring structure of the cohomology of a polyhedral product space
depends not only on all $i_k^*\colon H^*(X_k)\to H^*(A_k)$,
but also the character coproducts of $(X_k,A_k)$ defined in Theorem 2.8.
As an application, we give two polyhedral product spaces
${\cal Z}(K;X_1,A_1)$ and ${\cal Z}(K;X_2,A_2)$ in Example 7.13 such that
$A_1\simeq A_2$, $X_1\simeq X_2$ and the two cohomology homomorphisms
induced by inclusion are the same,
but the cohomology rings $H^*({\cal Z}(K;X_1,A_1))$ and $H^*({\cal Z}(K;X_2,A_2))$ are not isomorphic.
%\vspace{3mm}

\section{Character Coproduct}%\vspace{3mm}

\hspace*{5.5mm}
{\bf Notations and Conventions} All CW-complexes and CW-pairs in this paper have a base point.
For a CW-complex $X$, $(C_*(X),d)$ denotes
the cell chain complex. So $(C_*(X{\times}Y),d)\cong(C_*(X){\otimes}C_*(Y),d)$.
We use CW-complexes only to make the technical discussions easier.
\vspace{3mm}

{\bf Definition 2.1} A CW-pair $(X,A)$ is {\it homology split}
if $H_*(A)$, $H_*(X)$ and ${\rm ker}\,i_*$ are all free groups, where
$i_*:H_*(A)\!\to \!H_*(X)$ is induced by inclusion.
\vspace{3mm}

{\bf Definition 2.2} For a homology split $(X,A)$, the {\it character chain complex} $(C_*^\XX\!(X|A),d)$
is defined as follows.\vspace{-1mm}
$$C_*^\XX\!(X|A)\cong{\rm ker}\,i_*\oplus {\rm coker}\,i_*\oplus
{\rm im}\,i_*\oplus \Sigma\,{\rm ker}\,i_*,\vspace{-1mm}$$
where $\Sigma$ means uplifting the degree by $1$.
The restriction of $d$ on ${\rm ker}\,i_*\oplus {\rm coker}\,i_*\oplus
{\rm im}\,i_*$ is $0$ and the restriction of $d$ on $\Sigma\,{\rm ker}\,i_*$
is the desuspension isomorphism from $\Sigma\,{\rm ker}\,i_*$ to ${\rm ker}\,i_*$.
\vspace{3mm}

{\bf Theorem 2.3}\, {\it For a homology split $(X,A)$,
there is a quotient chain homotopy equivalence
$q\colon(C_*(X),d)\to(C_*^\XX\!(X|A),d)$
satisfying the following
commutative diagram
$$\begin{array}{ccc}
(C_*(A),d)&\stackrel{q'}{\longrightarrow}&H_*(A)\,{\cong}\,{\rm ker}\,i_*{\oplus}{\rm im}\,i_*\,\vspace{1mm}\\
\downarrow&&\downarrow\\
(C_*(X),d)&\stackrel{q}{\longrightarrow}&
(C_*^\XX\!(X|A),d),
  \end{array}$$
where $q'$ is the restriction of $q$ that is also a chain homotopy equivalence and
the two vertical homomorphisms are inclusions.
\vspace{2mm}

Proof}\, Take a representative $a_i$ in $C_*(A)$ for every generator of ${\rm ker}\,i_*$
and let $\overline a _i\in C_*(X)$ be any element such that $d\overline a _i= a_i$.
Take a representative $b_j$ in $C_*(A)$ for every generator of ${\rm im}\,i_*$.
Take a representative $c_k$ in $C_*(X)$ for every generator of ${\rm coker}\,i_*$.
So we may regard $H_*(A)$ as the chain subcomplex of $C_*(A)$ freely generated by all $a_i$'s and $b_j$'s
and regard $(C_*^\XX\!(X|A),d)$ as the chain subcomplex of $C_*(X)$ freely generated by all
$a_i$'s, $\overline a _i$'s, $b_j$'s and $c_k$'s.
Since all the homology groups are free, there
are free chain subcomplexes $F_*(A)$ of $C_*(A)$ and $F_*(X)$
of $C_*(X)$ such that\\
\hspace*{5mm}$(C_*(A),d)=(F_*(A){\oplus}H_*(A),d)$,\,\,\, $(C_*(X),d)=(F_*(X){\oplus}C_*^\XX\!(X|A),d)$\\
and $F_*(A)$ is a chain subcomplex of $F_*(X)$. Define $q,q'$ to be the projections.
Then we have the commutative diagram of the theorem.
Since $H_*(C_*(X))\cong H_*(X)\cong H_*(C_*^\XX\!(X|A))$, we have $H_*(F_*(X))\cong0$.
So $q$ is a chain homotopy equivalence.
Similarly, $q'$ is a chain homotopy equivalence.
\hfill$\Box$\vspace{3mm}

{\bf Definition 2.4} For a simplicial complex $K$ on $[m]$
and a sequence of CW-pairs $(\underline{X},\underline{A})=\{(X_k,A_k)\}_{k=1}^m$,
the {\it polyhedral product space} ${\cal Z}(K;\underline{X},\underline{A})$
is the subspace of $X_1{\times}{\cdots}{\times}X_m$ defined as follows.
For a subset $\tau$ of $[m]$, define
$$D(\tau)= Y_1{\times}\cdots{\times}Y_m,\quad Y_k=\left\{\begin{array}{cl}
X_k&{\rm if}\,\,k\in \tau, \\
A_k&{\rm if}\,\,k\not\in \tau.
\end{array}
\right.$$
Then ${\cal Z}(K;\underline{X},\underline{A})=\cup_{\tau\in K}D(\tau)$.
Denote by\vspace{1mm}\\
\hspace*{25mm}$i_k\colon H_*(A_k)\to H_*(X_k),\quad i_k^*\colon H^*(X_k)\to H^*(A_k)$\vspace{1mm}\\
the (co)homology homomorphisms induced by the inclusion maps.

If every $(X_k,A_k)$ is $(X,A)$,
then ${\cal Z}(K;\underline{X},\underline{A})$ is denoted by ${\cal Z}(K;X,A)$.
A ghost vertex $\{i\}\not\in K$ is allowed.
So ${\cal Z}(\{\emptyset\};\underline{X},\underline{A})= A_1{\times}{\cdots}{\times}A_m$.
The void complex $\{\}$ is inevitable and we define
${\cal Z}(\{\};\underline{X},\underline{A})=\emptyset$.
\vspace{3mm}

{\bf Definition 2.5} A polyhedral product space ${\cal Z}(K;\underline{X},\underline{A})$
is {\it homology split} if every pair $(X_k,A_k)$ is homology split.

The {\it character chain complex} $(C_*^{\XX_m}(K;\underline{X},\underline{A}),d)$
of the homology split ${\cal Z}(K;\underline{X},\underline{A})$
is the chain subcomplex of
$(C_*^\XX\!(X_1|A_1)\otimes\cdots\otimes C_*^\XX\!(X_m|A_m),d)$ defined as follows.
For a subset $\tau$ of $[m]$, define\vspace{-1mm}
$$(H_*(\tau),d)=(H_1{\otimes}\cdots{\otimes}H_m,d),\quad (H_k,d)=\left\{\begin{array}{cl}
(C_*^\XX\!(X_k|A_k),d)&{\rm if}\,\,k\in \tau,\vspace{1mm}\\
H_*(A_k)&{\rm if}\,\,k\not\in \tau.
\end{array}
\right.\vspace{-1mm}$$
Then $(C_*^{\XX_m}(K;\underline{X},\underline{A}),d)=(+_{\tau\in K}\,H_*(\tau),d)$.
\vspace{3mm}

{\bf Theorem 2.6} {\it For a homology split ${\cal Z}(K;\underline{X},\underline{A})$,
there is a quotient chain homotopy equivalence\,
$q_{(K;\underline{X},\underline{A})}\colon(C_*({\cal Z}(K;\underline{X},\underline{A})),d)
\to(C_*^{\XX_m}(K;\underline{X},\underline{A}),d)$.
\vspace{2mm}

Proof}\, Let $q_k,q'_k$ be as in Theorem 2.3.
For $\tau\subset[m]$, define\vspace{-2mm}
$$q_\tau= p_1{\otimes}{\cdots}{\otimes}p_m,\quad
p_k=\left\{\begin{array}{cl}
q_k&{\rm if}\,\,k\in\tau,\\
q'_k&{\rm if}\,\,k\notin\tau.
\end{array}\right.$$
Since all the homology groups are free, by K\"{u}nneth theorem,
$q_\tau$ is a chain homotopy equivalence from $C_*(D(\tau))$ to $H_*(\tau)$. So\vspace{-2mm}
$$q_{(K;\underline{X},\underline{A})}\!=\!+_{\tau\in{\scriptscriptstyle K}}\,q_\tau\colon C_*(\cup_{\tau\in K}D(\tau))=+_{\tau\in K}\,C_*(D(\tau))
\to+_{\tau\in K}\,H_*(\tau)\vspace{-2mm}$$
is a chain homotopy equivalence.
\hfill$\Box$\vspace{3mm}

The homology of $C_*^{\XX_m}(K;\underline{X},\underline{A})$ will be computed from the point of view of
diagonal tensor product in Section 4. The ring structure of
the cohomology of the polyhedral product spaces depends on the character coproduct
defined in the following two theorems.
\vspace{3mm}

{\bf Theorem 2.7}\, {\it  Let $(X,A)$ be a homology split CW-pair and\vspace{-2mm}
$$f_*\colon (C_*(X),d)\to(C_*(X{\times}X),d)\cong(C_*(X){\otimes}C_*(X),d)\vspace{-2mm}$$
be induced by a cellular map $f$ that is homotopic to the diagonal map of $X$.
Suppose the restriction of $f$ on $A$ is homotopic to the diagonal map of $A$.

There is a chain homomorphism\vspace{-2mm}
$$f_1\colon (C_*^\XX\!(X|A),d)\to(C_*^\XX\!(X|A){\otimes}C_*^\XX\!(X|A),d)\vspace{-2mm}$$
satisfying the following commutative diagram ($q$ as in Theorem 2.3)\vspace{-2mm}
$$\begin{array}{ccc}
(C_*(X),d)&\stackrel{f_*}{-\!\!\!-\!\!\!-\!\!\!\longrightarrow}& (C_*(X){\otimes}C_*(X),d)\,\vspace{1mm}\\
{\scriptstyle q}\downarrow\quad&&{\scriptstyle q{\otimes}q}\downarrow\quad\quad\quad\\
(C_*^\XX\!(X|A),d)&\stackrel{f_1}{-\!\!\!-\!\!\!-\!\!\!\longrightarrow}&
(C_*^\XX\!(X|A){\otimes}C_*^\XX\!(X|A),d),\,
  \end{array}\vspace{-1mm}$$
such that the restriction of $f_1$ on $H_*(A)$ is $\psi_A$,
the coproduct of $H_*(A)$ induced by the diagonal map of $A$.
\vspace{2mm}

Proof}\, Let everything be as in the proof of Theorem 2.3 and $C_*^\XX=C_*^\XX\!(X|A)$.
Then
$C_*(X){\otimes}C_*(X)=F_*(X){\otimes}C_*(X)\oplus C_*^\XX{\otimes}F_*(X)\oplus\, C_*^\XX{\otimes}C_*^\XX$
and $(q{\otimes}q)(F_*(X){\otimes}C_*(X){\oplus}C_*^\XX{\otimes}F_*(X))=0$.
Since $f_*$ is a chain homomorphism, we have $f_*(F_*(X))\subset F_*(X){\otimes}C_*(X){\oplus}C_*^\XX{\otimes}F_*(X)$.
So there is $f_1$ (depending on the choice of $F_*(X)$) satisfying the condition of the theorem.

Similarly, we have the commutative diagram of restrictions
$$\begin{array}{ccc}
(C_*(A),d)&\stackrel{f_*|_{C_*(A)}}{-\!\!\!-\!\!\!-\!\!\!\longrightarrow}& (C_*(A){\otimes}C_*(A),d)\,\vspace{1mm}\\
{\scriptstyle q'}\downarrow\quad&&{\scriptstyle q'{\otimes}q'}\downarrow\quad\quad\quad\,\\
H_*(A)&\stackrel{f_1|_{H_*(A)}}{-\!\!\!-\!\!\!-\!\!\!\longrightarrow}&
H_*(A){\otimes}H_*(A).\quad
  \end{array}$$
By definition, $f_1|_{H_*(A)}=\psi_A$.
\hfill$\Box$\vspace{3mm}

{\bf Theorem 2.8} {\it Let $(X,A)$ be a homology split CW-pair.
There is a unique chain homomorphism\vspace{-2mm}
$$\psi_{(X|A)}\colon(C_*^\XX\!(X|A),d)\to(C_*^\XX\!(X|A){\otimes}C_*^\XX\!(X|A),d)\vspace{-2mm}$$
satisfying the following three conditions.

i) \,$\psi_{(X|A)}$ is chain homotopic to all $f_1$ in Theorem 2.7.

ii) The restriction of $\psi_{(X|A)}$ on $H_*(A)\subset C_*^\XX\!(X|A)$ is $\psi_A$.

iii) Denote by $\alpha={\rm coker}\,i_*$, $\beta=\Sigma\,{\rm ker}\,i_*$,
$\gamma={\rm ker}\,i_*$, $\eta={\rm im}\,i_*$.
Then $\psi_{(X|A)}$ satisfies the following four conditions.

(1) $\psi_{(X|A)}(\eta)\subset \eta{\otimes}\eta\oplus\gamma{\otimes}\eta\oplus \eta{\otimes}\gamma
\oplus \gamma{\otimes}\gamma$.

(2) $\psi_{(X|A)}(\gamma)\subset \gamma{\otimes}\gamma\oplus \gamma{\otimes}\eta\oplus \eta{\otimes}\gamma$.

(3) $\psi_{(X|A)}(\beta)\subset\big(\beta{\otimes}\gamma\oplus \beta{\otimes}\eta\oplus\eta{\otimes}\beta\big)
\oplus\big(\alpha{\otimes}\alpha\oplus \alpha{\otimes}\eta\oplus
\eta{\otimes}\alpha\oplus \eta{\otimes}\eta\big)$.

(4) $\psi_{(X|A)}(\alpha)\subset \alpha{\otimes}\alpha\oplus \alpha{\otimes}\eta\oplus \eta{\otimes}\alpha
\oplus \eta{\otimes}\eta$.

$\psi_{(X|A)}$ is called the character coproduct of $(X,A)$.
\vspace{2mm}

Proof}\, Let $U\oplus(\alpha{\oplus}\eta){\otimes}(\alpha{\oplus}\eta)= C_*^\XX\!(X|A){\otimes}C_*^\XX\!(X|A)$.
By K\"{u}nneth theorem, $H_*(U)=0$ and $H_*(C_*^\XX\!(X|A){\otimes}C_*^\XX\!(X|A))=(\alpha{\oplus}\eta){\otimes}(\alpha{\oplus}\eta)$.
For $f_1$ in Theorem 2.7, construct a chain homomorphism\\ \hspace*{25mm}$\psi\colon(C_*^\XX\!(X|A),d)\to(C_*^\XX\!(X|A){\otimes}C_*^\XX\!(X|A),d)$\\
and a chain homotopy\\
\hspace*{25mm}$s\colon(C_*^\XX\!(X|A),d)\to(\Sigma C_*^\XX\!(X|A){\otimes}C_*^\XX\!(X|A),d)$\\
such that $ds+sd= f_1{-}\psi$ as follows.

For $x\in\gamma{\oplus}\eta= H_*(A)$, define $\psi(x)= f_1(x)$ and $s(x)=0$.
Then $(ds+sd)(x)=(f_1{-}\psi)(x)$ and $\psi$ naturally satisfies (1) and (2) on $\gamma{\oplus}\eta$.

For a generator $b\in\beta$ with $db=c$ and $f_1(c)=\Sigma\,c'_i{\otimes}c''_i$,
suppose $f_1(b)= x+y$,
where $x\in U$, $y\in(\alpha{\oplus}\eta){\otimes}(\alpha{\oplus}\eta)$.
Define $\psi(b)=\Sigma\,b'_i{\otimes}b''_i+y$, where $db'_i= c'_i$, $b''_i= c''_i$ if $c'_i\in\gamma$ and
$b'_i= c'_i$, $db''_i=(-1)^{|b'_i|}c''_i$ if $c'_i\notin\gamma$.
Then $f_1(b){-}\psi(b)\in U$ and $d(f_1(b){-}\psi(b))= f_1(c)-\psi(c)=0$.
Since $H_*(U)=0$, there is $z\in U$ such that $dz= f_1(b)-\psi(b)$.
Define $s(b)= z$. Then $(ds+sd)(b)=(f_1{-}\psi)(b)$ and $\psi$ satisfies (3) on $\beta$.

For a generator $a\in\alpha$, $d(f_1(a))=0$ implies that $f_1(a)= u+v$, where
$u\in U$, $du=0$ and $v\in(\alpha{\oplus}\eta){\otimes}(\alpha{\oplus}\eta)$.
Since $H_*(U)=0$, there is $w$ such that $dw= u$.
Define $\psi(a)= v$ and $s(a)= w$. Then $(ds+sd)(a)=(f_1{-}\psi)(a)$ and
$\psi$ satisfies (4) on $\alpha$ .

Suppose $f_1$ in Theorem 2.8 is replaced by another $f'_1$ and $\psi'$ is constructed as above for $f'_1$. By i), $\psi|_{\gamma\oplus\eta}=\psi'|_{\gamma{\oplus}\eta}$.
By definition, if $\psi(b)\neq\psi'(b)$ for $b\in\beta$, then $\psi(b){-}\psi'(b)\in (\alpha{\oplus}\eta)\otimes(\alpha{\oplus}\eta)
= H_*(C_*^\XX\!(X|A){\otimes}C_*^\XX\!(X|A))$.
This contradicts $\psi\simeq\psi'$. So $\psi|_{\beta}=\psi'|_{\beta}$.
Similarly, $\psi|_{\alpha}=\psi'|_{\alpha}$.
Thus, $\psi=\psi'$. This implies that $\psi$ does not depend on a choice of $f_1$.
So we may denote the unique $\psi_{(X|A)}$ by $\psi$.
\hfill$\Box$\vspace{3mm}

{\bf Example 2.9} Consider the map\vspace{-2mm}
$$f\colon S^3\stackrel{\mu'}{-\!\!\!-\!\!\!\longrightarrow}S^3\vee S^3
\stackrel{g\,\vee\, 1}{-\!\!\!-\!\!\!\longrightarrow} S^2\vee S^3,\vspace{-2mm}$$
where $\mu'$ is the coproduct of the co-$H$-space $S^3$, $g\colon S^3\to S^2$ is the Hopf bundle and $1$ is the identity map of $S^3$.
Let $A_1= A_2= S^2\vee S^3$, $X_1= S^2\vee CS^3$, $X_2= C_f$, where $C$ means the cone of
a space and $C_f$ means the mapping cone of $f$.
By definition, $(C_*^\XX\!(X_1|A_1),d)=(C_*^\XX\!(X_2|A_2),d)$ and they are freely generated
by $1,a,b,\overline b$ with $d\overline b= b$, $|1|=0$, $|a|=2$, $|b|=3$.
So the equality $H_*(X_i)\cong H_*(C_*^\XX\!(X_i|A_i))\cong H_*(S^2)$ implies $X_1\simeq X_2\simeq S^2$.
By definition,
$\psi_{(X_1|A_1)}(\overline b)=1{\otimes}\overline b+\overline b{\otimes}1$,
$\psi_{(X_2|A_2)}(\overline b)=1{\otimes}\overline b+\overline b{\otimes}1+a{\otimes}a$.
So $\psi_{(X_1|A_1)}\not\simeq\psi_{(X_2|A_2)}$.
As we will see in Example 7.13, $H^*({\cal Z}(K;X_1,A_1))$ and $H^*({\cal Z}(K;X_2,A_2))$ are isomorphic as groups
but not isomorphic as rings.
\vspace{3mm}

In the following definition, the tensor product of coproducts are defined as follows.
For $\psi_1(a)=\Sigma_i\,a'_i{\otimes}a''_i$ and
$\psi_2(b)=\Sigma_j\,b'_j{\otimes}b''_j$,
$(\psi_1{\otimes}\psi_2)(a{\otimes}b)=\Sigma_{i,j}\,
(-1)^{|a''_i||b'_j|}(a'_i{\otimes}b'_j){\otimes}(a''_i{\otimes}b''_j)$.\vspace{2mm}

{\bf Definition 2.10} For a homology split space ${\cal Z}(K;\underline{X},\underline{A})$,
the {\it character coproduct} is a chain homomorphism\vspace{-2mm}
$$\psi_{(K;\underline{X},\underline{A})}
\colon(C_*^{\XX_m}(K;\underline{X},\underline{A}),d)\to(C_*^{\XX_m}(K;\underline{X},\underline{A})
{\otimes}C_*^{\XX_m}(K;\underline{X},\underline{A}),d)\vspace{-2mm}$$
defined as follows. For $\tau\subset[m]$, define ($H_*(\tau)$ as in Definition 2.5)\vspace{-2mm}
$$\psi_\tau\colon (H_*(\tau),d)\to(H_*(\tau){\otimes}H_*(\tau),d)\vspace{-2mm}$$
by\vspace{-2mm}
$$\psi_\tau=\psi_1{\otimes}{\cdots}{\otimes}\psi_m,\quad
\psi_k=\left\{\begin{array}{cl}
\psi_{(X_k|A_k)}&{\rm if}\,\,k\in\tau,\\
\psi_{A_k}&{\rm if}\,\,k\notin\tau.
\end{array}\right.$$
Then $\psi_{(K;\underline{X},\underline{A})}=+_{\tau\in K}\,\psi_{\tau}$.
\vspace{3mm}

{\bf Theorem 2.11} {\it Let ${\cal Z}(K;\underline{X},\underline{A})$ be homology split.
Then for a cellular map $f\colon{\cal Z}(K;\underline{X},\underline{A})\to
{\cal Z}(K;\underline{X},\underline{A}){\times}{\cal Z}(K;\underline{X},\underline{A})$
that is homotopic to the diagonal map,
the following diagram is homotopy commutative
$$\begin{array}{ccc}
(C_*({\cal Z}(K;\underline{X},\underline{A})),d)&\stackrel{f_*}{\longrightarrow}
&(C_*({\cal Z}(K;\underline{X},\underline{A})){\otimes}
C_*({\cal Z}(K;\underline{X},\underline{A})),d)\,\,\,\,\,\,\,\vspace{1mm}\\
{\scriptstyle q_{(K;\underline{X},\underline{A})}}\downarrow\hspace{13mm}&&{\scriptstyle q_{(K;\underline{X},\underline{A})}{\otimes}q_{(K;\underline{X},\underline{A})}}\downarrow\hspace{35mm}\,\\
(C_*^{\XX_m}(K;\underline{X},\underline{A}),d)
&\stackrel{\psi_{(K;\underline{X},\underline{A})}}{-\!\!\!-\!\!\!-\!\!\!\longrightarrow}
&(C_*^{\XX_m}(K;\underline{X},\underline{A}){\otimes}
C_*^{\XX_m}(K;\underline{X},\underline{A}),d).\quad
\end{array}\vspace{1mm}$$

Proof}\, Let $f_k$ be a cellular map homotopic to the diagonal map of $X_k$ such that the restriction of $f_k$
on $A_k$ is homotopic to the diagonal map of $A_k$. We may take $f$ to be the restriction
of $f_1{\times}{\cdots}{\times}f_m$.
Let $f_\tau\colon (C_*(D(\tau)),d)\to(C_*(D(\tau){\times}D(\tau)),d)$ be the restriction of $f_*$
on $C_*(D(\tau))$.
Then we have the following commutative diagram\vspace{-2mm}
$$\begin{array}{ccc}
(C_*(D(\tau)),d)&\stackrel{f_\tau}{\longrightarrow}&(C_*(D(\tau){\times}D(\tau)),d)\,\\
{\scriptstyle q_\tau}\downarrow\quad&&{\scriptstyle q_\tau{\otimes}q_\tau}\downarrow\hspace{10mm}\\
(H_*(\tau),d)&\stackrel{g_\tau}{\longrightarrow}&(H_*(\tau){\otimes}H_*(\tau)),d),
\end{array}
$$
where
$$g_\tau= g_1{\otimes}{\cdots}{\otimes}g_m,\quad
g_k=\left\{\begin{array}{cl}
(f_k)_1&{\rm if}\,\,k\in\tau,\vspace{2mm}\\
\psi_{A_k}&{\rm if}\,\,k\notin\tau,
\end{array}\right.$$
and $(f_k)_1$ corresponds to $f_1$ in Theorem 2.7 for $(X,A)=(X_k,A_k)$.
Define $g_{(K;\underline{X},\underline{A})}=+_{\tau\in{\scriptscriptstyle K}}\,g_\tau$.
So we have the following commutative diagram
$$\begin{array}{ccc}
(C_*({\cal Z}(K;\underline{X},\underline{A})),d)&\stackrel{f_*}{\longrightarrow}
&(C_*({\cal Z}(K;\underline{X},\underline{A})){\otimes}
C_*({\cal Z}(K;\underline{X},\underline{A})),d)\,\,\,\,\,\,\,\vspace{1mm}\\
{\scriptstyle q_{(K;\underline{X},\underline{A})}}\downarrow\hspace{13mm}&&{\scriptstyle q_{(K;\underline{X},\underline{A})}{\otimes}q_{(K;\underline{X},\underline{A})}}\downarrow\hspace{35mm}\,\\
(C_*^{\XX_m}(K;\underline{X},\underline{A}),d)
&\stackrel{g_{(K;\underline{X},\underline{A})}}{-\!\!\!-\!\!\!-\!\!\!\longrightarrow}
&(C_*^{\XX_m}(K;\underline{X},\underline{A}){\otimes}
C_*^{\XX_m}(K;\underline{X},\underline{A}),d).\quad
\end{array}\vspace{2mm}$$
By Theorem 2.8, $g_\tau\simeq \psi_\tau$.
So $g_{(K;\underline{X},\underline{A})}\simeq\psi_{(K;\underline{X},\underline{A})}$.
\hfill$\Box$

\section{Diagonal Tensor Product}

\hspace{5.5mm}{\bf Definition 3.1} A {\it group $A_*^\Lambda$ indexed by the set} $\Lambda$
is a direct sum over $\Lambda$ of graded groups $A_*^\Lambda=\oplus_{\alpha\in\Lambda}\,A_*^{\alpha}$.

A {\it chain complex $(C_*^\Lambda,d)$ indexed by the set} $\Lambda$
is a direct sum over $\Lambda$ of chain complexes $(C_*^\Lambda,d)=\oplus_{\alpha\in\Lambda}(C_*^{\alpha},d)$.

A {\it cochain complex $(C^*_\Lambda,\delta)$ indexed by the set} $\Lambda$
is a  direct sum over $\Lambda$ of cochain complexes $(C^*_\Lambda,\delta)=\oplus_{\alpha\in\Lambda}(C^*_{\alpha},\delta)$.
\vspace{3mm}

{\bf Lemma 3.2} {\it Every subgroup or quotient group of a group indexed by $\Lambda$
is naturally a group indexed by $\Lambda$.

For $A_*^\Lambda=\oplus_{\alpha\in\Lambda}\,A_*^\alpha$ with $\Lambda$ a finite set,
its dual $A^*_\Lambda={\rm Hom}(A_*^\Lambda,\zz)$ is a group indexed by $\Lambda$ such that $A^*_\Lambda=\oplus_{\alpha\in\Lambda}\,A^*_\alpha$ with $A^*_\alpha={\rm Hom}(A_*^\alpha,\zz)$.

The conclusion also holds for chain and cochain complexes.
\vspace{2mm}

Proof}\, Suppose $B$ is a subgroup of $A_*^\Lambda=\oplus_{\alpha\in\Lambda}\,A_*^{\alpha}$,
then $B=\oplus_{\alpha\in\Lambda}\,B{\cap}A_*^{\alpha}$ is indexed by $\Lambda$.
$A_*^\Lambda/B=\oplus_{\alpha\in\Lambda}\,A_*^{\alpha}/(B{\cap}A_*^{\alpha})$ is indexed by $\Lambda$.
\hfill$\Box$\vspace{3mm}

{\bf Definition 3.3} Let $A_*^\Lambda=\oplus_{\alpha\in\Lambda}\,A_*^\alpha$ and $B_*^\Gamma=\oplus_{\beta\in\Gamma}\,B_*^\beta$ be groups indexed by $\Lambda$ and $\Gamma$ respectively.
Their {\it tensor product group}\vspace{1mm}\\
\hspace*{40mm}$A_*^\Lambda\otimes B_*^\Gamma=\oplus_{(\alpha,\beta)\in\Lambda\times\Gamma}\,A_*^{\alpha}\otimes B_*^{\beta}\vspace{1mm}$\\
is naturally indexed by $\Lambda\times\Gamma$.

Let $A_*^\Lambda=\oplus_{\alpha\in\Lambda}\,A_*^\alpha$ and
$B_*^\Lambda=\oplus_{\alpha\in\Lambda}\,B_*^\alpha$ be groups indexed by the same set $\Lambda$.
Their {\it diagonal tensor product} (with respect to $\Lambda$) is the group indexed by $\Lambda$
given by\\
\hspace*{32mm}$A_*^\Lambda\otimes_\Lambda B_*^\Lambda=\oplus_{\alpha\in\Lambda}\,C_*^\alpha,\,\,\,\,
C_*^\alpha= A_*^{\alpha}\otimes B_*^{\alpha}.$

The {\it tensor product} and {\it diagonal tensor product} of indexed (co)chain complexes
are defined by replacing all the groups in the above
definitions by (co)chain complexes.
\vspace{3mm}

We have to deal with tensor product and diagonal tensor product simultaneously. For example,
$(A_*^\Lambda{\otimes}_\Lambda B_*^\Lambda)\otimes({A'}_*^\Lambda{\otimes}_\Lambda {B'}_*^\Lambda)$.
So we must have the following convention.
\vspace{3mm}

{\bf Convention} For $A_*^\Lambda$ and $B_*^\Lambda$ indexed by the same set,
we use $a\widehat\otimes b$ to denote the element of $A_*^\Lambda{\otimes}_\Lambda B_*^\Lambda$
and $a{\otimes}b$ to denote the element of $A_*^\Lambda{\otimes}B_*^\Lambda$.
Precisely, for $a\in A_*^\alpha\subset A_*^\Lambda$ and $b\in B_*^\beta\subset B_*^\Lambda$,
define $a{\widehat\otimes}b= a{\otimes}b\in A_*^\alpha{\otimes}B_*^\alpha\subset A_*^\Lambda{\otimes_\Lambda}B_*^\Lambda$
if $\alpha=\beta$ and $a{\widehat\otimes}b=0$ if $\alpha\neq\beta$.
\vspace{3mm}

{\bf Theorem 3.4} {\it There is an indexed group (or complex) isomorphism\vspace{1mm}\\
\hspace*{35mm}$(A_{\Lambda_1}{\otimes}_{\Lambda_1}B_{\Lambda_1}){\otimes}{\cdots}{\otimes}
(A_{\Lambda_m}{\otimes}_{\Lambda_m}B_{\Lambda_m})$\vspace{1mm}\\
\hspace*{30.5mm}$\cong
(A_{\Lambda_1}{\otimes}{\cdots}{\otimes}A_{\Lambda_m}){\otimes}_{\Lambda_1\times\cdots\times\Lambda_m}
(B_{\Lambda_1}{\otimes}{\cdots}{\otimes}B_{\Lambda_m})$,\vspace{1.5mm}\\
where $(-)_{\Lambda_i}$ means $(-)_*^{\Lambda_i}$ or $(-)^*_{\Lambda_i}$.
\vspace{2mm}

Proof}\, The restriction of the factor-permuting isomorphism\vspace{-2mm}
$$\phi\colon A_{\Lambda_1}{\otimes}B_{\Lambda_1}{\otimes}{\cdots}{\otimes}
A_{\Lambda_m}{\otimes}B_{\Lambda_m}\stackrel{\cong}{\longrightarrow}
A_{\Lambda_1}{\otimes}{\cdots}{\otimes}A_{\Lambda_m}{\otimes}
B_{\Lambda_1}{\otimes}{\cdots}{\otimes}B_{\Lambda_m}\vspace{-2mm}$$
on the subgroup $(A_{\Lambda_1}{\otimes}_{\Lambda_1}B_{\Lambda_1})\otimes\cdots\otimes
(A_{\Lambda_m}{\otimes}_{\Lambda_m}B_{\Lambda_m})$
is just the isomorphism of the theorem.
Precisely, for $a_i\in A_{\Lambda_i}$, $b_i\in B_{\Lambda_i}$, \vspace{-2mm}
$$\phi((a_1{\otimes}b_1){\otimes}{\cdots}{\otimes}(a_m{\otimes}b_m))=
(-1)^s(a_1{\otimes}{\cdots}{\otimes}a_m){\otimes}(b_1{\otimes}{\cdots}{\otimes}b_m),\vspace{-2mm}$$
where $s=\Sigma_{i=2}^m(|b_1|{+}{\cdots}{+}|b_{i-1}|)|a_i|$. Then\vspace{-2mm}
$$\hat\phi\colon (A_{\Lambda_1}{\otimes}_{\Lambda_1}B_{\Lambda_1}){\otimes}{\cdots}{\otimes}
(A_{\Lambda_m}{\otimes}_{\Lambda_m}B_{\Lambda_m})\hspace{40mm}\vspace{-2mm}$$
$$\hspace{20mm}\stackrel{\cong}{\longrightarrow}
(A_{\Lambda_1}{\otimes}{\cdots}{\otimes}A_{\Lambda_m}){\otimes}_{\Lambda_1\times\cdots\times\Lambda_m}
(B_{\Lambda_1}{\otimes}{\cdots}{\otimes}B_{\Lambda_m})\vspace{-2mm}$$
satisfies that for $a_i\in A_*^{\alpha_i}$, $b_i\in B_*^{\alpha_i}$, $\alpha_i\in\Lambda_i$,\vspace{-2mm}
$$\hspace*{9mm}\hat\phi((a_1{\widehat\otimes}b_1){\otimes}{\cdots}{\otimes}(a_m{\widehat\otimes}b_m))=
(-1)^s(a_1{\otimes}{\cdots}{\otimes}a_m){\widehat\otimes}(b_1{\otimes}{\cdots}{\otimes}b_m).\hspace{10mm}\Box$$

\section{Homology and Cohomology Group}

\hspace*{5.5mm}{\bf Notations and Conventions} In this paper, all definitions and theorems have a dual analogue
and all the dual proofs are omitted. The index set\vspace{-2mm}
$$\XX_m=\{(\sigma\!,\omega)\,|\,\sigma\!,\omega\subset[m],\,\sigma\!\cap\!\omega=\emptyset\}.\vspace{-2mm}$$
Define $\XX$ to be $\XX_1=\{(\emptyset,\emptyset),(\emptyset,\{1\}),(\{1\},\emptyset)\}$.
Then $\XX_m=\XX{\times}{\cdots}{\times}\XX$ ($m$ fold) by the following 1-1 correspondence\vspace{-2mm}
$$(\sigma,\omega)\to(s_1,{\cdots},s_m),\quad
s_k=\left\{\begin{array}{cl}
{\scriptstyle(\{1\},\emptyset)}&{\rm if}\,\,k\in\sigma,\\
{\scriptstyle(\emptyset,\{1\})}&{\rm if}\,\,k\in\omega,\\
{\scriptstyle(\,\emptyset,\,\emptyset\,)}&{\rm otherwise}.
\end{array}\right.\vspace{3mm}$$

{\bf Definition 4.1} Let $(X,A)$ be a homology split CW-pair.

The {\it homology group $H_*^\XX\!(X,A)$ indexed by $\XX$} is given by\vspace{-2mm}
$$H_*^{\emptyset,\emptyset}(X,A)={\rm im}\,i_*,\,\,
H_*^{\emptyset,\{1\}}(X,A)={\rm ker}\,i_*,\,\,
H_*^{\{1\},\emptyset}(X,A)={\rm coker}\,i_*.\vspace{-2mm}$$

The {\it character chain complex $(C_*^\XX\!(X|A),d)$} in Definition 2.2 is a chain complex indexed by $\XX$
given by\vspace{-2mm}
$$(C_*^{\emptyset,\emptyset}(X|A),d)={\rm im}\,i_*,\,\,
(C_*^{\{1\},\emptyset}(X|A),d)={\rm coker}\,i_*,\vspace{-2mm}$$
$$(C_*^{\emptyset,\{1\}}(X|A),d)=({\rm ker}\,i_*\oplus\Sigma\,{\rm ker}\,i_*,d).\vspace{-2mm}$$
It is obvious that $H_*^\XX\!(X,A)$ is a trivial chain subcomplex of $C_*^\XX\!(X|A)$.

Dually, by Lemma 3.2, the {\it cohomology group $H^*_\XX(X,A)$ indexed by $\XX$} is
the dual group of $H_*^\XX\!(X,A)$ given by\vspace{-2mm}
$$H^*_{\emptyset,\emptyset}(X,A)={\rm im}\,i^*,\,\,
H^*_{\emptyset,\{1\}}(X,A)={\rm coker}\,i^*,\,\,
H^*_{\{1\},\emptyset}(X,A)={\rm ker}\,i^*.\vspace{-2mm}$$

The {\it character cochain complex $(C^*_\XX(X|A),\delta)$} is the dual of $(C_*^\XX\!(X|A),d)$
given by\vspace{-2mm}
$$(C^*_{\emptyset,\emptyset}(X|A),\delta)={\rm im}\,i^*,\,\,
(C^*_{\{1\},\emptyset}(X|A),\delta)={\rm ker}\,i^*,\vspace{-2mm}$$
$$(C^*_{\emptyset,\{1\}}(X|A),\delta)=({\rm coker}\,i^*\oplus\Sigma{\rm coker}\,i^*,\delta).$$

{\bf Definition 4.2} Denote  by $\zz(x_1,{\cdots},x_n)$
the free abelian group generated by $x_1,{\cdots},x_n$.

$(T_*^\XX,d)$ is a chain complex indexed by $\XX$ defined as follows.\vspace{-2mm}
$$(T_*^{\emptyset,\emptyset},d)=\zz(\eta),\,\,(T_*^{\emptyset,\{1\}},d)=(\zz(\beta,\gamma),d),
(T_*^{\{1\},\emptyset},d)=\zz(\alpha),\vspace{-2mm}$$
where $|\alpha|=|\gamma|=|\eta|=0$, $|\beta|=1$, $d\beta=\gamma$.

$(T^*_\XX,\delta)$ is the dual cochain complex of $(T_*^\XX,d)$ defined as follows.\vspace{-2mm}
$$(T^*_{\emptyset,\emptyset},\delta)=\zz(\eta),\,\,(T^*_{\emptyset,\{1\}},\delta)=(\zz(\beta,\gamma),\delta),
(T^*_{\{1\},\emptyset},\delta)=\zz(\alpha),\vspace{-2mm}$$
where we use the same symbol to denote a generator and its dual generator. So $\delta\gamma=\beta$.
\vspace{3mm}

{\bf Theorem 4.3} {\it Let $(X,A)$ be a homology split CW-pair.

There is an isomorphism of chain complexes indexed by $\XX$\vspace{-2mm}
$$\phi\colon (C_*^\XX\!(X|A),d)\stackrel{\cong}{\longrightarrow}
(T_*^\XX\,{\otimes}_\XX\, H_*^\XX\!(X,A),d)\vspace{-2mm}$$
such that the restriction of $\phi$ on $H_*(A)$ is the following isomorphism
$$\phi'\colon H_*(A)\stackrel{\cong}{\longrightarrow}
S_*^\XX\,{\otimes}_\XX\, H_*^\XX\!(X,A),\vspace{-2mm}$$
where $S_*^\XX$ is the subgroup $\zz(\gamma,\eta)$ (by Lemma 3.2, $S_*^\XX$ is indexed by $\XX$).

Dually, the dual \vspace{-2mm}
$$\phi^*\colon(T^*_\XX\,{\otimes}_\XX\, H^*_\XX(X,A),\delta)\stackrel{\cong}{\longrightarrow} (C^*_\XX(X|A),\delta)\vspace{-2mm}$$
of $\phi$ is an isomorphism of cochain complexes indexed by $\XX$.
\vspace{2mm}

Proof}\, Define $\phi$ as follows.\vspace{3mm}\\
\hspace*{31mm}
\begin{tabular}{|c|c|c|c|c|}
\hline
{\rule[-2mm]{0mm}{6mm}$\quad x\in$}&${\rm coker}\,i_*$&$\Sigma\,{\rm ker}\,i_*$
&${\rm ker}\,i_*$&${\rm im}\,i_*$\\
\hline
{\rule[-2mm]{0mm}{7mm}$\phi(x)=$}&$\alpha\,\widehat\otimes\,x$&$\beta\,\widehat\otimes\,dx$
&$\gamma\,\widehat\otimes\,x$&$\eta\,\widehat\otimes\,x$\\
\hline
\end{tabular}
\vspace{3mm}\\
It is obvious that $\phi$ is a chain isomorphism.
\hfill$\Box$\vspace{3mm}

{\bf Definition 4.4} Let $(\underline{X},\underline{A})=\{(X_k,A_k)\}_{k=1}^m$ be a sequence of
CW-pairs such that every pair $(X_k,A_k)$ is homology split.

The {\it homology group $H_*^{\XX_m}(\underline{X},\underline{A})$ indexed by $\XX_m$} is given by\vspace{-2mm}
$$H_*^{\XX_m}(\underline{X},\underline{A})=
H_*^\XX\!(X_1,A_1)\otimes{\cdots}\otimes H_*^\XX\!(X_m,A_m).\vspace{-2mm}$$
Denote $H_*^{\XX_m}(\underline{X},\underline{A})=\oplus_{(\sigma\!,\,\omega)\in\XX_m}
H_*^{\sigma\!,\,\omega}(\underline{X},\underline{A})$. Then\vspace{-2mm}
$$H_*^{\sigma\!,\,\omega}(\underline{X},\underline{A})
= H_1{\otimes}{\cdots}{\otimes}H_m,\quad
H_k=\left\{\begin{array}{cl}
{\rm coker}\,i_k&{\rm if}\, k\in\sigma,\vspace{1mm}\\
{\rm ker}\,i_k&{\rm if}\, k\in\omega,\vspace{1mm}\\
{\rm im}\,i_k&{\rm otherwise}.
\end{array}\right.\vspace{-1mm}$$

The {\it character chain complex $(C_*^{\XX_m}(\underline{X}|\underline{A}),d)$
indexed by $\XX_m$} is given by\vspace{-2mm}
$$(C_*^{\XX_m}(\underline{X}|\underline{A}),d)=
(C_*^\XX\!(X_1|A_1)\otimes{\cdots}\otimes C_*^\XX\!(X_m|A_m),d).\vspace{-2mm}$$

Dually, the {\it cohomology group $H^*_{\XX_m}(\underline{X},\underline{A})$ indexed by $\XX_m$}
is given by\vspace{-2mm}
$$H^*_{\XX_m}(\underline{X},\underline{A})=
H^*_\XX(X_1,A_1)\otimes{\cdots}\otimes H^*_\XX(X_m,A_m).\vspace{-2mm}$$
Then $H^*_{\XX_m}(\underline{X},\underline{A})=\oplus_{(\sigma\!,\,\omega)\in\XX_m}
H^*_{\sigma\!,\,\omega}(\underline{X},\underline{A})$ with
$$H^*_{\sigma\!,\,\omega}(\underline{X},\underline{A})
= H^1{\otimes}{\cdots}{\otimes}H^m,\quad
H^k=\left\{\begin{array}{cl}
{\rm ker}\,i_k^*&{\rm if}\, k\in\sigma,\vspace{1mm}\\
{\rm coker}\,i_k^*&{\rm if}\, k\in\omega,\vspace{1mm}\\
{\rm im}\,i_k^*&{\rm otherwise}.
\end{array}\right.\vspace{3mm}$$
The {\it character cochain complex $(C^*_{\XX_m}(\underline{X}|\underline{A}),\delta)$
indexed by $\XX_m$} is the dual cochain complex
$(C^*_\XX(X_1|A_1)\otimes{\cdots}\otimes C^*_\XX(X_m|A_m),\delta)$
of $(C_*^{\XX_m}(\underline{X}|\underline{A}),d)$.
\vspace{3mm}

{\bf Definition 4.5} Let $K$ be a simplicial complex on $[m]$. Denote by\vspace{-2mm}
$$(T_*^{\XX_m},d)=(T_*^\XX{\otimes}{\cdots}{\otimes}T_*^\XX,d),\quad
(T^*_{\XX_m},\delta)=(T^*_\XX{\otimes}{\cdots}{\otimes}T^*_\XX,\delta)\,\,
(m\,{\rm fold}). \vspace{-2mm}$$

The {\it total chain complex $(T_*^{\XX_m}(K),d)$ of $K$ indexed by $\XX_m$} is
the chain subcomplex of $(T_*^{\XX_m},d)$ defined as follows.
For a subset $\tau$ of $[m]$, define
$$(T_*(\tau),d)=(T_1{\otimes}{\cdots}{\otimes}T_m,d),\quad
T_k=\left\{\begin{array}{cl}
T_*^\XX&{\rm if}\,\,k\in\tau,\vspace{1mm}\\
S_*^\XX&{\rm if}\,\,k\notin\tau.
\end{array}\right.$$
Then $(T_*^{\XX_m}(K),d)=(+_{\tau\in K}\,T_*(\tau),d)$.
The {\it total homology group of $K$ indexed by $\XX_m$} is defined to be
$H_*^{\XX_m}(K)= H_*(T_*^{\XX_m}(K))$.

Dually, the {\it total cochain complex $(T^*_{\XX_m}(K),\delta)$ of $K$ indexed by $\XX_m$}
is the dual of $(T_*^{\XX_m}(K),d)$.
The {\it total cohomology group of $K$ indexed by $\XX_m$} is defined to be $H^*_{\XX_m}(K)= H^*(T^*_{\XX_m}(K))$.
\vspace{3mm}

Since the character chain complex $(C_*^{\XX_m}(K;\underline{X},\underline{A}),d)$
in Definition 2.5 is a chain subcomplex of $(C_*^{\XX_m}(\underline{X}|\underline{A}),d)$,
it is a chain complex indexed by $\XX_m$ by Lemma 3.2.
Dually, the character cochain complex $(C^*_{\XX_m}(K;\underline{X},\underline{A}),\delta)$
as the dual of $(C_*^{\XX_m}(K;\underline{X},\underline{A}),d)$ is indexed by $\XX_m$.
\vspace{3mm}

{\bf Theorem 4.6} {\it For a homology split ${\cal Z}(K;\underline{X},\underline{A})$,
there is an isomorphism of chain complexes indexed by $\XX_m$\vspace{-2mm}
$$\phi_{(K;\underline{X},\underline{A})}\colon(C_*^{\XX_m}(K;\underline{X},\underline{A}),d)
\stackrel{\cong}{\longrightarrow}
(T_*^{\XX_m}(K){\otimes}_{\XX_m}H_*^{\XX_m}(\underline{X},\underline{A}),d)\vspace{-2mm}$$
such that the dual\vspace{-2mm}
$$\phi_{(K;\underline{X},\underline{A})}^*\colon
(T^*_{\XX_m}(K){\otimes}_{\XX_m}C^*_{\XX_m}(\underline{X},\underline{A}),\delta)
\stackrel{\cong}{\longrightarrow}(C^*_{\XX_m}(K;\underline{X},\underline{A}),\delta)\vspace{-2mm}$$
is an isomorphism of cochain complexes indexed by $\XX_m$.
\vspace{2mm}

Proof}\, Let\, $\phi_k\colon C_*^\XX\!(X_k|A_k)\to T_*^\XX\otimes_\XX H_*^\XX\!(X_k,A_k)$\,
and\, its\, restriction\, $\phi'_k\colon$ $H_*(A_k)\to S_*^\XX{\otimes}_\XX H_*^\XX\!(X_k,A_k)$
be as in Theorem 4.3. For $\tau\subset[m]$, define\vspace{-2mm} $$\phi_\tau=\lambda_1{\otimes}{\cdots}{\otimes}\lambda_m,\quad
\lambda_k=\left\{\begin{array}{cl}
\phi_k&{\rm if}\,\,k\in\tau,\\
\phi'_k&{\rm if}\,\,k\notin\tau.
\end{array}\right.$$
Then $\phi_\tau\colon H_*(\tau)\to
\big(T_1{\otimes}_\XX H_*^\XX\!(X_1,A_1)\big){\otimes}{\cdots}{\otimes}
\big(T_m{\otimes}_\XX H_*^\XX\!(X_m,A_m)\big)$ is an isomorphism. By Theorem 3.4,\vspace{-2mm}
$$\big(T_1{\otimes}_\XX H_*^\XX\!(X_1,A_1)\big){\otimes}{\cdots}{\otimes}
\big(T_m{\otimes}_\XX H_*^\XX\!(X_m,A_m)\big)
\cong T_*(\tau){\otimes}_{\XX_m}H_*^{\XX_m}(\underline{X},\underline{A}).\vspace{-1mm}$$
Identifying the two chain complexes, we have an isomorphism\vspace{-1mm}
$$\phi_\tau\colon(H_*(\tau),d)\stackrel{\cong}{\longrightarrow}
\big(T_*(\tau){\otimes}_{\XX_m} H_*^{\XX_m}(\underline{X},\underline{A}),d\big).\vspace{-1mm}$$
Specifically, for $\tau=[m]$, we have the isomorphism\vspace{-2mm}
$$\phi_{[m]}\colon(C_*^{\XX_m}(\underline{X}|\underline{A}),d)\stackrel{\cong}{\longrightarrow}
\big(T_*^{\XX_m}{\otimes}_{\XX_m} H_*^{\XX_m}(\underline{X},\underline{A}),d\big).\vspace{-1mm}$$
Define $\phi_{(K;\underline{X},\underline{A})}{=}+_{\tau\in K}\phi_\tau\colon+_{\tau\in K}H_*(\tau)
\to(+_{\tau\in K}\,T_*(\tau)){\otimes}_{\XX_m}H_*^{\XX_m}(\underline{X},\underline{A})$.
As a restriction of $\psi_{[m]}$ into its image,
$\phi_{(K;\underline{X},\underline{A})}$ is an isomorphism.
\hfill$\Box$\vspace{3mm}

{\bf Theorem 4.7} {\it  Let $K$ be a simplicial complex on $[m]$. Denote the
groups indexed by $\XX_m$ in Definition 4.5 by\vspace{-2mm}
$$T_*^{\XX_m}(K)=\oplus_{(\sigma,\omega)\in\XX_m}\,T_*^{\sigma\!,\,\omega}(K),\quad
H_*^{\XX_m}(K)=\oplus_{(\sigma,\omega)\in\XX_m}\,H_*^{\sigma\!,\,\omega}(K),\vspace{-2mm}$$
$$T^*_{\XX_m}(K)=\oplus_{(\sigma,\omega)\in\XX_m}\,T^*_{\sigma\!,\,\omega}(K),\quad
H^*_{\XX_m}(K)=\oplus_{(\sigma,\omega)\in\XX_m}\,H^*_{\sigma\!,\,\omega}(K).$$
Then for every $(\sigma\!,\,\omega)\in\XX_m$,\vspace{-2mm}
$$(T_*^{\sigma\!,\,\omega}(K),d)\cong(\Sigma\w C_*(K_{\sigma\!,\,\omega}),d),\quad
H_*^{\sigma\!,\,\omega}(K)\cong\w H_{*-1}(K_{\sigma\!,\,\omega}),\vspace{-2mm}$$
$$(T^*_{\sigma\!,\,\omega}(K),\delta)\cong(\Sigma\w C^*(K_{\sigma\!,\,\omega}),\delta),\quad
H^*_{\sigma\!,\,\omega}(K)\cong\w H^{*-1}(K_{\sigma\!,\,\omega}),$$
where $\Sigma\w C$ means the augmented simplicial (co)chain complex with degree uplifted by $1$,
$K_{\sigma\!,\,\omega}=({\rm link}_{_K}\sigma)|_\omega=\{\tau\,\,|\,\,\tau\subset\omega,\,\tau{\cup}\sigma\in K\}$
if $\sigma\in K$ and $K_{\sigma\!,\,\omega}=\{\}$ (the void complex) if $\sigma\notin K$.
\vspace{2mm}

Proof}\, Denote by $t_{A,B,C,D}$ the generator $t_1{\otimes}{\cdots}{\otimes}t_m$ of $T_*^{\XX_m}$
such that\vspace{1mm}\\
\hspace*{7mm}$A=\{k\,|\,t_k\!=\!\alpha\}$, $B=\{k\,|\,t_k\!=\!\beta\}$, $C=\{k\,|\,t_k\!=\!\gamma\}$, $D=\{k\,|\,t_k\!=\!\eta\}$.\vspace{1mm}\\
Then for $\tau\subset[m]$,
$T(\tau)=\zz(\{t_{A,B,C,D}\}_{A\cup B\subset\,\tau})$
and $T_*^{\XX_m}=\oplus_{(\sigma,\omega)\in\XX_m}\,T_*^{\sigma,\omega}$ with
$T_*^{\sigma,\omega}=\zz(\{t_{A,B,C,D}\}
_{A=\sigma,\,B\cup C=\omega})$. So\vspace{1mm}\\
\hspace*{14mm}$T_*^{\sigma\!,\,\omega}(K)=+_{\tau\in K}\,T_*^{\sigma\!,\,\omega}\cap T_*(\tau)
=\zz(\{t_{\sigma,B,C,D}\}
_{\sigma\cup B\,\in\, K,\,\,B\cup C\,=\,\omega})$.\vspace{1mm}\\
The 1-1 correspondence $B\to
t_{\sigma,\,B,\,\omega{\setminus}B,\,[m]{\setminus}(\sigma{\cup}\omega)}$
for all $B\in K_{\sigma\!,\,\omega}$ induces a chain isomorphism
from $(\Sigma\w C_*(K_{\sigma\!,\,\omega}),d)$ to $(T_*^{\sigma\!,\,\omega}(K),d)$ if $\sigma\in K$
and $T_*^{\sigma,\omega}(K)=\Sigma\w C_*(\{\})=0$ if $\sigma\notin K$.
\hfill$\Box$\vspace{3mm}

{\bf Theorem 4.8} {\it For a homology split $M={\cal Z}(K;\underline{X},\underline{A})$,\vspace{-1mm}
$$H_*(M)\cong
H_*^{\XX_m}(K){\otimes}_{\XX_m}H_*^{\XX_m}(\underline{X},\underline{A})
\cong\oplus_{(\sigma,\omega)\in\XX_m}\,
\w H_{*-1}(K_{\sigma,\omega}){\otimes}H_*^{\sigma\!,\,\omega}
(\underline{X},\underline{A}),\vspace{-1mm}$$
$$H^*(M)\cong
H^*_{\XX_m}(K){\otimes}_{\XX_m}H^*_{\XX_m}(\underline{X},\underline{A})
\cong\oplus_{(\sigma,\omega)\in\XX_m}\,
\w H^{*-1}(K_{\sigma,\omega}){\otimes}H^*_{\sigma\!,\,\omega}
(\underline{X},\underline{A}).\vspace{-1mm}$$

The conclusion holds for all polyhedral product spaces if the (co)homology group is taken over a field.
}\vspace{2mm}

{\it Proof}\, We have $H_*(M)\cong
H_*(T_*^{\XX_m}(K){\otimes}_{\XX_m}H_*^{\XX_m}(\underline{X},\underline{A}))$
by Theorem 2.6 and Theorem 4.6.
Since
$H_*^{\XX_m}(\underline{X},\underline{A})$ is a free and trivial chain complex, we have
$H_*(T_*^{\XX_m}(K){\otimes}_{\XX_m}H_*^{\XX_m}(\underline{X},\underline{A}))
\cong H_*^{\XX_m}(K)
{\otimes}_{\XX_m}H_*^{\XX_m}(\underline{X},\underline{A})$.

All definitions and proofs in this paper have natural generalizations to (co)homology
over a field and the condition that every $(X_k,A_k)$ is homology split is superfluous in this case.
\hfill$\Box$\vspace{3mm}

{\bf Example 4.9} We compute the (co)homology group of ${\cal Z}(K;S^r,S^p)$, where $p<r$ and $S^p$ is a subcomplex of $S^r$
by any tame embedding. Then\vspace{-2mm}
$${\rm im}\,i_*\cong\zz,\quad\quad{\rm coker}\,i_*\cong\zz,\quad\quad{\rm ker}\,i_*\cong\zz.\vspace{-2mm}$$
So $H_*^{\sigma\!,\,\omega}(\underline{S^r},\underline{S^p})\cong\zz$ for all $(\sigma\!,\,\omega)\in\XX_m$.
Identifying
$H_*^{\sigma\!,\,\omega}(K){\otimes}H_*^{\sigma\!,\,\omega}(\underline{S^r},\underline{S^p})$ with
$H_*^{\sigma\!,\,\omega}(K)$ (degree uplifted), we have\vspace{-2mm}
$$H_*({\cal Z}(K;S^r,S^p))\cong\oplus_{(\sigma,\omega)\in\XX_m}\,\w H_{*-r|\sigma|-p|\omega|-1}(K_{\sigma\!,\,\omega}),\vspace{-2mm}$$
$$H^*({\cal Z}(K;S^r,S^p))\cong\oplus_{(\sigma\!,\,\omega)\in\XX_m}\,\w H^{*-r|\sigma|-p|\omega|-1}(K_{\sigma\!,\,\omega}).$$

\section{Diagonal Tensor Product of Algebras and Coalgebras}

\hspace*{5.5mm}{\bf Definition 5.1} A {\it coalgebra $(A_*^\Lambda,\psi)$ indexed by $\Lambda$}
is an indexed group $A_*^\Lambda$ with a coproduct $\psi\colon A_*^\Lambda\to A_*^\Lambda{\otimes}A_*^\Lambda$ that is a group homomorphism (may not keep degree or be coassociative).

A {\it subcoalgebra} $(B_*^\Lambda,\psi')$ of $(A_*^\Lambda,\psi)$ is a coalgebra such that
$B_*^\Lambda$ is a subgroup of $A_*^\Lambda$ and $\psi'$ is the restriction of $\psi$.

Dually, \,an \,{\it algebra\, $(A^*_\Lambda,\pi)$\, indexed\, by\, $\Lambda$\,} is\, an\, indexed\, group\,
$A^*_\Lambda$\, with\, a\, product\,
$\pi\colon A^*_\Lambda{\otimes}A^*_\Lambda\to A^*_\Lambda$ that is a group homomorphism
(may not keep degree or be associative).

An ideal $I^*_\Lambda$ of $A^*_\Lambda$ is a subgroup such that
$\pi(I^*_\Lambda{\otimes}A^*_\Lambda+A^*_\Lambda{\otimes}I^*_\Lambda)\subset I^*_\Lambda$.
$\pi$ induces a product $\pi'\colon A^*_\Lambda/I^*_\Lambda\otimes A^*_\Lambda/I^*_\Lambda\to A^*_\Lambda/I^*_\Lambda$
and $(A^*_\Lambda/I^*_\Lambda,\pi')$ is called a {\it quotient algebra} of $(A^*_\Lambda,\psi)$.
\vspace{3mm}

{\bf Remark} We have to deal with a coalgebra $(A_*^\Lambda,\psi)$ that is also a chain complex $(A_*^\Lambda,d)$.
In this paper, the coproduct $\psi$ in this case is a chain homomorphism from
$(A_*^\Lambda,d)$ to $(A_*^\Lambda{\otimes}A_*^\Lambda,d)$
and all subcoalgebras $(B_*^\Lambda,\psi')$ of $(A_*^\Lambda,\psi)$ are assumed to be chain subcomplexes.
Analogue conventions hold for algebras that is also a cochain complex.
\vspace{3mm}

{\bf Lemma 5.2} {\it Let $(A_*^\Lambda,\psi)$ be a coalgebra such that $\Lambda$ is finite and
$A_*^\Lambda$ is free. Then the dual group $A^*_\Lambda={\rm Hom}(A_*^\Lambda,\zz)$ is an algebra with
product $\pi\colon A^*_\Lambda{\otimes}A^*_\Lambda\to A^*_\Lambda$ the dual of $\psi$
define by $\pi(f{\otimes}g)(a)=(f{\otimes}g)(\psi(a))$ for all $a\in A_*^\Lambda$
and $f,g\in A^*_\Lambda$.
$(A^*_\Lambda,\pi)$ is called the dual algebra of $(A_*^\Lambda,\psi)$.
For a subcoalgebra $(B_*^\Lambda,\psi')$ of $(A_*^\Lambda,\psi)$,
the dual algebra $(B^*_\Lambda,\pi')$ of $(B_*^\Lambda,\psi')$
is a quotient algebra of $(A^*_\Lambda,\pi)$.
\vspace{2mm}

Proof}\,
For a subcoalgebra $(B_*^\Lambda,\psi')$, let
$I^*_\Lambda=\{f\in A^*_\Lambda\,|\, f(b)=0\,\,{\rm for\,\,all}\,\,b\in B_*^\Lambda\}$.
Then $I^*_\Lambda$ is an ideal of $A^*_\Lambda$ such that ${\rm Hom}(B_*^\Lambda,\zz)= A^*_\Lambda/I^*_\Lambda$.
\hfill$\Box$\vspace{3mm}

{\bf Definition 5.3} For coalgebra $(A_*^\Lambda\!=\!\oplus_{\alpha\in\Lambda}\,A_*^\alpha,\psi)$,
the coproduct $\psi$ is determined by all its {\it restriction coproduct} defined as follows.
For $a\in A_*^\alpha$ and every $\beta,\gamma\in\Lambda$, there is a unique
$b_{\beta,\gamma}\in A_*^\beta{\otimes}A_*^\gamma$ such that $\psi(a)=\Sigma_{\beta,\gamma\in\Lambda}b_{\beta,\gamma}$.
The correspondence $a\to b_{\beta,\gamma}$ is the group homomorphism\vspace{-2mm}
$$\psi^\alpha_{\beta,\gamma}\colon A^\alpha_*\stackrel{i}{\to}A_*^\Lambda
\stackrel{\psi}{\longrightarrow}A_*^\Lambda{\otimes}A_*^\Lambda\stackrel{p}{\to}
A^\beta_*{\otimes}A^\gamma_*,\vspace{-2mm}$$
where $i$ is the inclusion and $p$ is the projection. Every $\psi^\alpha_{\beta,\gamma}$
is called a restriction coproduct of $\psi$.
$\psi$ is defined if and only if all its restriction coproducts are defined.

Dually, for algebra $(A^*_\Lambda\!=\!\oplus_{\alpha\in\Lambda}\,A^*_\alpha,\pi)$,
the product $\pi$ is determined by all its {\it restriction product} defined as follows.
For $b\in A^*_\beta$, $c\in A^*_\gamma$ and every $\alpha\in\Lambda$, there is a unique
$a_\alpha\in A^*_\alpha$ such that
$\pi(b{\otimes}c)=\Sigma_{\alpha\in\Lambda}\,a_\alpha$.
The correspondence $b{\otimes}c\to a_\alpha$ is the group homomorphism\vspace{-2mm}
$$\pi_\alpha^{\beta,\gamma}\colon A_\beta^*{\otimes}A_\gamma^*\stackrel{i}{\to}
A^*_\Lambda{\otimes}A^*_\Lambda\stackrel{\pi}{\longrightarrow}A^*_\Lambda\stackrel{p}{\to} A_\alpha^*,\vspace{-2mm}$$
where $i$ is the inclusion and $p$ is the projection. Every \,$\pi_\alpha^{\beta,\gamma}$\,
is\, called\, a\, restriction\, product\, of $\pi$.
$\pi$ is defined if and only if all its restriction products are defined.
\vspace{3mm}

{\bf Definition 5.4} Let $(A_*^\Lambda,\psi_1)$ and $(B_*^\Gamma,\psi_2)$
be two coalgebras. Their {\it tensor product coalgebra} $(A_*^\Lambda{\otimes} B_*^\Gamma,\psi_1{\otimes}\psi_2)$
is defined as follows. Suppose for $\alpha\in\Lambda$ and $a\in A_*^\alpha$, $\psi_1(a)=\Sigma_i\,a'_i{\otimes}a''_i$
with every $a'_i{\otimes}a''_i\in A_*^{\alpha'}{\otimes}A_*^{\alpha''}$ for some $\alpha'\!,\alpha''\in\Lambda$.
Suppose for $\beta\in\Gamma$ and $b\in B_*^\beta$, $\psi_2(b)=\Sigma_j\,b'_j{\otimes}b''_j$
with every $b'_j{\otimes}b''_j\in B_*^{\beta'}{\otimes}B_*^{\beta''}$ for some $\beta'\!,\beta''\in\Gamma$.
Define\vspace{1mm}\\
\hspace*{24mm}$(\psi_1{\otimes}\psi_2)(a{\otimes}b)=\Sigma_{i,j}\,
(-1)^{|a''_i||b'_j|}(a'_i{\otimes}b'_j){\otimes}(a''_i{\otimes}b''_j)$.\vspace{1mm}\\
Equivalently, every restriction coproduct of $\psi_1{\otimes}\psi_2$ satisfies\vspace{1mm}\\
\hspace*{30mm}$(\psi_1{\otimes}\psi_2)^{(\alpha,\beta)}_{(\alpha'\!,\beta'),(\alpha''\!,\beta'')}=
(\psi_1)^{\alpha}_{\alpha'\!,\alpha''}\otimes(\psi_2)^{\beta}_{\beta'\!,\beta''}.$\vspace{1mm}

Dually, let $(A^*_\Lambda,\pi_1)$ and $(B^*_\Gamma,\pi_2)$ be two algebras. Their {\it tensor product
algebra} $(A^*_\Lambda{\otimes} B^*_\Gamma,\pi_1{\otimes}\pi_2)$ is
defined as follows. Suppose for $\alpha'\!\in A^*_{\alpha'}$ and  $a''\in A^*_{\alpha''}$,
$\pi_1(a'{\otimes}a'')=\Sigma_i\,a_i$  with every $a_i\in A^*_\alpha$ for some $\alpha\in\Lambda$.
Suppose for $\beta'\in B^*_{\beta'}$ and $b''\in B^*_{\beta''}$, $\pi_2(b'{\otimes}b'')=\Sigma_j\,b_j$
with every $b_j\in B^*_\beta$ for some $\beta\in\Gamma$. Define \vspace{1mm}\\
\hspace*{23mm}$(\pi_1{\otimes}\pi_2)((a'{\otimes}b'){\otimes}(a''{\otimes}b''))
=(-1)^{|a''||b'|}(\Sigma_{i,j}\,a_i{\otimes}b_j)$.\vspace{1mm}\\
Equivalently, every restriction product of $\pi_1{\otimes}\pi_2$ satisfies\vspace{1mm}\\
\hspace*{32mm}$(\pi_1{\otimes}\pi_2)_{(\alpha,\beta)}^{(\alpha'\!,\beta'),(\alpha''\!,\beta'')}=
(\pi_1)_{\alpha}^{\alpha'\!,\alpha''}\otimes
(\pi_2)_{\beta}^{\beta'\!,\beta''}$.
\vspace{3mm}

{\bf Definition 5.5} Let $(A_*^\Lambda,\psi_1)$ and $(B_*^\Lambda,\psi_2)$
be two coalgebras indexed by the same set. Their {\it diagonal tensor product coalgebra}
$(A_*^\Lambda{\otimes}_\Lambda B_*^\Lambda,\psi_1{\otimes}_\Lambda\psi_2)$
is defined as follows. Suppose for $a\in A_*^\alpha$ and $b\in B_*^\alpha$, $\psi_1(a)=\Sigma_i\,a'_i{\otimes}a''_i$
with every $a'_i{\otimes}a''_i\in A_*^{\alpha'}{\otimes}A_*^{\alpha''}$ for some $\alpha'\!,\alpha''\in\Lambda$
and $\psi_2(b)=\Sigma_j\,b'_j{\otimes}b''_j$
with every $b'_j{\otimes}b''_j\in B_*^{\beta'}{\otimes}B_*^{\beta''}$ for some $\beta'\!,\beta''\in\Lambda$.
Define\vspace{1mm}\\
\hspace*{24mm}$(\psi_1{\otimes}_\Lambda\psi_2)(a{\widehat\otimes}b)=\Sigma_{i,j}\,
(-1)^{|a''_i||b'_j|}(a'_i{\widehat\otimes}b'_j){\otimes}(a''_i{\widehat\otimes}b''_j)$.\vspace{1mm}\\
Equivalently, every restriction coproduct of $\psi_1{\otimes}_\Lambda\psi_2$ satisfies\vspace{1mm}\\
\hspace*{34mm}$(\psi_1{\otimes}_\Lambda\psi_2)^{\alpha}_{\alpha'\!,\alpha''}=
(\psi_1)^{\alpha}_{\alpha'\!,\alpha''}\otimes(\psi_2)^{\alpha}_{\alpha'\!,\alpha''}.$\vspace{1mm}

Dually, let $(A^*_\Lambda,\pi_1)$ and $(B^*_\Lambda,\pi_2)$ be two algebras indexed by the same set.
Their {\it diagonal tensor product algebra} $(A^*_\Lambda{\otimes}_\Lambda B^*_\Lambda,\pi_1{\otimes}_\Lambda\pi_2)$ is
defined as follows. Suppose for $\alpha'\!,\alpha''\in\Lambda$ and  $a'{\otimes}a''\in A^*_{\alpha'}{\otimes}A^*_{\alpha''}$
and  $b'{\otimes}b''\in B^*_{\alpha'}{\otimes}B^*_{\alpha''}$,
$\pi_1(a'{\otimes}a'')=\Sigma_i\,a_i$  with every $a_i\in A^*_\alpha$ for some $\alpha\in\Lambda$
and $\pi_2(b'{\otimes}b'')=\Sigma_j\,b_j$
with every $b_j\in B^*_\beta$ for some $\beta\in\Lambda$. Define \vspace{1mm}\\
\hspace*{23mm}$(\pi_1{\otimes}_\Lambda\pi_2)((a'{\widehat\otimes}b'){\otimes}(a''{\widehat\otimes}b''))
=(-1)^{|a''||b'|}(\Sigma_{i,j}\,a_i{\widehat\otimes}b_j)$.\vspace{1mm}\\
Equivalently, every restriction product of $\pi_1{\otimes}_\Lambda\pi_2$ satisfies\vspace{1mm}\\
\hspace*{34mm}$(\pi_1{\otimes}_\Lambda\pi_2)_{\alpha}^{\alpha'\!,\alpha''}=
(\pi_1)_{\alpha}^{\alpha'\!,\alpha''}\otimes
(\pi_2)_{\alpha}^{\alpha'\!,\alpha''}$.
\vspace{3mm}

The properties of diagonal tensor product are very different from that of tensor product.
For example, for non-associative algebras $(A^*_\Lambda,\pi_1)$ and $(B^*_\Lambda,\pi_2)$,
$(A^*_\Lambda{\otimes}_\Lambda B^*_\Lambda,\pi_1{\otimes}_\Lambda\pi_2)$ may be an associative, commutative algebra.
For $(A^*_\Lambda,\pi_1)\not\cong(A'^*_\Lambda,\pi'_1)$ and  $(B^*_\Lambda,\pi_2)\not\cong(B'^*_\Lambda,\pi'_2)$,
there might be an isomorphism $(A^*_\Lambda{\otimes}_\Lambda B^*_\Lambda,\pi_1{\otimes}_\Lambda\pi_2)\cong
(A'^*_\Lambda{\otimes}_\Lambda B'^*_\Lambda,\pi'_1{\otimes}_\Lambda\pi'_2)$.
This is because for $a{\otimes}b\neq 0$ in a tensor product group,
$a{\widehat\otimes}b$ may be $0$ in the diagonal tensor product group.
\vspace{3mm}

{\bf Theorem 5.6} {\it There is a (co)algebra isomorphism\vspace{-2mm}
$$\Big((A_{\Lambda_1}{\otimes}_{\Lambda_1}B_{\Lambda_1}){\otimes}{\cdots}{\otimes}
(A_{\Lambda_m}{\otimes}_{\Lambda_m}B_{\Lambda_m}),(\varphi_1{\otimes}_{\Lambda_1}\varphi'_1)
{\otimes}{\cdots}{\otimes}(\varphi_m{\otimes}_{\Lambda_m}\varphi'_m)\Big)\vspace{-2mm}$$
$$\cong
\Big((A_{\Lambda_1}{\otimes}{\cdots}{\otimes}A_{\Lambda_m}){\otimes}_{\Lambda}
(B_{\Lambda_1}{\otimes}{\cdots}{\otimes}B_{\Lambda_m}),(\varphi_1{\otimes}{\cdots}{\otimes}\varphi_m)
{\otimes}_{\Lambda}
(\varphi'_1{\otimes}{\cdots}{\otimes}\varphi'_m)
\Big),\vspace{-2mm}$$
where $\Lambda=\Lambda_1{\times}{\cdots}{\times}\Lambda_m$ and $((-)_{\Lambda_i},\varphi_i)$ means indexed
(co)algebra.\vspace{2mm}

Proof}\, The group isomorphism in Theorem 3.4 is naturally a (co)algebra isomorphism.
\hfill$\Box$

\section{Cohomology Algebra}

\hspace*{5.5mm} In Section 4, we proved
$H^*(C^*_{\XX_m}(K;\underline{X},\underline{A}))\cong
H^*_{\XX_m}(K){\otimes}_{\XX_m}H^*_{\XX_m}(\underline{X},\underline{A})$.
In this section, we will prove that this group isomorphism is a ring isomorphism.
We define coproducts $\psi$ and $\psi_{(X,A)}$ on $T_*^\XX$ and $H_*^\XX\!(X,A)$ respectively
and get a coalgebra isomorphism\\
\hspace*{19mm}$(C_*^\XX\!(X|A),\psi_{(X|A)})\cong(T_*^\XX{\otimes}_\XX H_*^\XX\!(X,A),\,\psi{\otimes}_\XX\psi_{(X,A)})$\\
that induces the cohomology ring isomorphism.
\vspace{3mm}

{\bf Definition 6.1} The {\it universal coproduct}
$\psi\colon(T_*^\XX,d)\to(T_*^\XX{\otimes}T_*^\XX,d)$ is defined as follows.

$\psi(\eta)=\eta{\otimes}\eta+\gamma{\otimes}\eta+\eta{\otimes}\gamma+\gamma{\otimes}\gamma$.

$\psi(\gamma)=\gamma{\otimes}\gamma+\gamma{\otimes}\eta+\eta{\otimes}\gamma$.

$\psi(\beta)=\beta{\otimes}\gamma+\beta{\otimes}\eta+\eta{\otimes}\beta+
\alpha{\otimes}\alpha+\alpha{\otimes}\eta+\eta{\otimes}\alpha+\eta{\otimes}\eta$.

$\psi(\alpha)=\alpha{\otimes}\alpha+\alpha{\otimes}\eta+\eta{\otimes}\alpha+\eta{\otimes}\eta$.

The {\it normal coproduct}
$\w\psi\colon(T_*^\XX,d)\to(T_*^\XX{\otimes}T_*^\XX,d)$ is defined as follows.

$\w\psi(\eta)=\eta{\otimes}\eta+\eta{\otimes}\gamma+\gamma{\otimes}\eta+\gamma{\otimes}\gamma$.

$\w\psi(\gamma)=\gamma{\otimes}\gamma+\gamma{\otimes}\eta+\eta{\otimes}\gamma$.

$\w\psi(\beta)=\beta{\otimes}\gamma+\beta{\otimes}\eta+\eta{\otimes}\beta$.

$\w\psi(\alpha)=\alpha{\otimes}\alpha+\alpha{\otimes}\eta+\eta{\otimes}\alpha+\eta{\otimes}\eta$.

The {\it special coproduct}
$\overline\psi\colon(T_*^\XX,d)\to(T_*^\XX{\otimes}T_*^\XX,d)$ is defined as follows.

$\overline\psi(\eta)=\eta{\otimes}\eta$.

$\overline\psi(\gamma)=\gamma{\otimes}\eta+\eta{\otimes}\gamma$.

$\overline\psi(\beta)=\beta{\otimes}\eta+\eta{\otimes}\beta$.

$\overline\psi(\alpha)=\alpha{\otimes}\eta+\eta{\otimes}\alpha$.

By Lemma 5.2, the dual of the above coalgebras are algebras and are respectively denoted by
$(T^*_\XX,\pi),(T^*_\XX,\w\pi),(T^*_\XX,\overline\pi)$.
\vspace{3mm}

{\bf Definition 6.2} For a homology split pair $(X,A)$, the coproduct\vspace{-2mm}
$$\psi_{(X,A)}\colon H_*^\XX\!(X,A)\to H_*^\XX\!(X,A){\otimes}H_*^\XX\!(X,A)\vspace{-2mm}$$
is defined as follows, \,where  $\alpha={\rm coker}\,i_*$, $\beta=\Sigma\,{\rm ker}\,i_*$,
$\gamma={\rm ker}\,i_*$, $\eta={\rm im}\,i_*$ as in Theorem 2.8.

(1) $\psi_{(X,A)}(x)=\psi_{(X|A)}(x)$ for all $x\in \alpha\oplus\eta$.

(2) For a generator $x\in\gamma$, there is a unique generator $\overline x\in\beta$ such that $d\,\overline x= x$.
Suppose $\psi_{(X|A)}(\overline x)= z+y$ with $z\in\beta{\otimes}\gamma\oplus \beta{\otimes}\eta\oplus\eta{\otimes}\beta$
and $y\in\alpha{\otimes}\alpha\oplus \alpha{\otimes}\eta\oplus\eta{\otimes}\alpha\oplus \eta{\otimes}\eta$.
Then define $\psi_{(X,A)}(x)=\psi_{(X|A)}(x)+y$.

The {\it homology coalgebra indexed by $\XX$} of $(X,A)$ is $(H_*^\XX\!(X,A),\psi_{(X,A)})$.

Dually, by Lemma 5.2, the {\it cohomology algebra indexed by $\XX$} of $(X,A)$
is the dual algebra $(H^*_\XX(X,A),\pi_{(X,A)})$.
\vspace{3mm}

{\bf Definition 6.3} Let  $(X,A)$ be a homology split pair.

$(X,A)$ is called {\it normal} if the character coproduct satisfies

$\psi_{(X|A)}(\eta)\subset
\eta{\otimes}\eta\oplus\eta{\otimes}\gamma\oplus\gamma{\otimes}\eta\oplus\gamma{\otimes}\gamma$,

$\psi_{(X|A)}(\gamma)\subset\gamma{\otimes}\gamma\oplus\gamma{\otimes}\eta\oplus\eta{\otimes}\gamma$,

$\psi_{(X|A)}(\beta)\subset\beta{\otimes}\gamma\oplus\beta{\otimes}\eta\oplus\eta{\otimes}\beta$,

$\psi_{(X|A)}(\alpha)\subset\alpha{\otimes}\alpha\oplus\alpha{\otimes}\eta\oplus\eta{\otimes}\alpha
\oplus\eta{\otimes}\eta$.

$(X,A)$ is called {\it special} if the character coproduct satisfies

$\psi_{(X|A)}(\eta)\subset\eta{\otimes}\eta$,

$\psi_{(X|A)}(\gamma)\subset\gamma{\otimes}\eta\oplus\eta{\otimes}\gamma$,

$\psi_{(X|A)}(\beta)\subset\beta{\otimes}\eta\oplus\eta{\otimes}\beta$,

$\psi_{(X|A)}(\alpha)\subset\alpha{\otimes}\eta\oplus\eta{\otimes}\alpha$.
\vspace{3mm}

{\bf Theorem 6.4} {\it The group isomorphisms $\phi$ and $\phi^*$ in Theorem 4.3 is an (co)algebra isomorphism, i.e.,\vspace{-2mm}
$$(C_*^\XX\!(X|A),\psi_{(X|A)})\cong(T_*^\XX{\otimes}_\XX H_*^\XX\!(X,A),\,\psi{\otimes}_\XX\psi
_{(X,A)}),\vspace{-2mm}$$
$$(C^*_\XX(X|A),\pi_{(X|A)})\cong(T^*_\XX{\otimes}_\XX H^*_\XX(X,A),\,\pi{\otimes}_\XX\pi
_{(X,A)}).$$

If $(X,A)$ is normal (or special), then $\psi$ and $\pi$ can be replaced by
$\w\psi$ and $\w\pi$ (or $\overline\psi$ and $\overline\pi$) respectively.
\vspace{2mm}

Proof}\, We use the following symbols to denote elements of the corresponding groups.
\vspace{2mm}\\
\hspace*{5mm}
\begin{tabular}{|c|c|c|c|c|}
\hline
{\rule[-2mm]{0mm}{6mm}elements\,\,of}&${\rm coker}\,i_*$&$\Sigma\,{\rm ker}\,i_*$
&${\rm ker}\,i_*$&${\rm im}\,i_*$\\
\hline
{\rule[-2mm]{0mm}{7mm}symbols}&$a,a'_1,a''_1,\cdots$&$b,b'_1,b''_1,\cdots$
&$c,c'_1,c''_1,\cdots$&$e,e'_1,e''_1,\cdots$\\
\hline
\end{tabular}
\vspace{3mm}\\
All $\Sigma$ are omitted, i.e.,
$\Sigma\, x{\otimes}y$ is denoted by $x{\otimes}y$.
For $a,b,c,e\in C_*^\XX\!(X|A)$ such that $db=c$, $db'_i=c'_i$, $db''_i= c''_i$, suppose

$\psi_{(X|A)}(e)= e'_1{\otimes}e''_1+e'_2{\otimes}c''_2+c'_3{\otimes}e''_3+c'_4{\otimes}c''_4$,

$\psi_{(X|A)}(c)= c'_5{\otimes}c''_5+c'_6{\otimes}e''_6+e'_7{\otimes}c''_7$,

$\psi_{(X|A)}(b)= b'_5{\otimes}c''_5+b'_6{\otimes}e''_6+(-1)^{|e'_7|}e'_7{\otimes}b''_7\,+\,
a'_8{\otimes}a''_8\,+\,a'_9{\otimes}e''_9\,+\,e'_{10}{\otimes}a''_{10}\,+\,e'_{11}{\otimes}e''_{11}$,

$\psi_{(X|A)}(a)= a'_{12}{\otimes}a''_{12}+a'_{13}{\otimes}e''_{13}+e'_{14}{\otimes}a''_{14}+e'_{15}{\otimes}e''_{15}$.\\
Then

$\psi_{(X,A)}(e)= e'_1{\otimes}e''_1+e'_2{\otimes}c''_2+c'_3{\otimes}e''_3+c'_4{\otimes}c''_4$,

$\psi_{(X,A)}(c)= c'_5{\otimes}c''_5+c'_6{\otimes}e''_6+e'_7{\otimes}c''_7+\,
a'_8{\otimes}a''_8\,+\,a'_9{\otimes}e''_9\,+\,e'_{10}{\otimes}a''_{10}\,+\,e'_{11}{\otimes}e''_{11}$,

$\psi_{(X,A)}(a)= a'_{12}{\otimes}a''_{12}+a'_{13}{\otimes}e''_{13}+e'_{14}{\otimes}a''_{14}+e'_{15}{\otimes}e''_{15}$.

For simplicity, $x{\otimes}y$ is abbreviated to $xy$ and $x\widehat\otimes y$ is abbreviated to
$x{\scriptstyle\wedge}y$ in the following computation.

\hspace*{20mm}$(\psi{\otimes}_\XX\psi_{(X,A)})(\phi(e))=(\psi{\otimes}_\XX\psi_{(X,A)})(\eta{\scriptstyle\wedge} e)\\
\hspace*{20mm}=(\eta\eta{+}\gamma\eta{+}\eta\gamma{+}\gamma\gamma){\scriptstyle\wedge}
(e'_1e''_1{+}e'_2c''_2{+}c'_3e''_3{+}c'_4c''_4)\\
\hspace*{20mm}=(\eta{\scriptstyle\wedge} e'_1)(\eta{\scriptstyle\wedge} e''_1)
{+}(\eta{\scriptstyle\wedge} e'_2)(\gamma{\scriptstyle\wedge} c''_2)
{+}(\gamma{\scriptstyle\wedge} c'_3)(\eta{\scriptstyle\wedge} e''_3)
{+}(\gamma{\scriptstyle\wedge} c'_4)(\gamma{\scriptstyle\wedge} c''_4)\\
\hspace*{20mm}=\phi(e'_1){\otimes}\phi(e''_1)
{+}\phi(e'_2){\otimes}\phi(c''_2){+}\phi(c'_3){\otimes}\phi(e''_3)
{+}\phi(c'_4){\otimes}\phi(c''_4)\\
\hspace*{20mm}=(\phi{\otimes}\phi)(\psi_{(X|A)}(e))$,

\hspace*{20mm}$(\psi{\otimes}_\XX\psi_{(X,A)})(\phi(c))=(\psi{\otimes}_\XX\psi_{(X,A)})(\gamma{\scriptstyle\wedge} c)\\
\hspace*{20mm}=(\gamma\gamma{+}\gamma\eta{+}\eta\gamma){\scriptstyle\wedge}
(c'_5c''_5{+}c'_6e''_6{+}e'_7c''_7{+}a'_8a''_8{+}a'_9e''_9{+}e'_{10}a''_{10}{+}e'_{11}e''_{11})\\
\hspace*{20mm}=(\gamma{\scriptstyle\wedge} c'_5)(\gamma{\scriptstyle\wedge} c''_5){+}
(\gamma{\scriptstyle\wedge} c'_6)(\eta{\scriptstyle\wedge} e''_6){+}
(\eta{\scriptstyle\wedge} e'_7)(\gamma{\scriptstyle\wedge} c''_7)\\
\hspace*{20mm}=\phi(c'_5){\otimes}\phi(c''_5){+}\phi(c'_6){\otimes}\phi(e''_6){+}
\phi(e'_7){\otimes}\phi(c''_7)\\
\hspace*{20mm}=(\phi{\otimes}\phi)(\psi_{(X|A)}(c))$,

\hspace*{20mm}$(\psi{\otimes}_\XX\psi_{(X,A)})(\phi(b))=(\psi{\otimes}_\XX\psi_{(X,A)})(\beta{\scriptstyle\wedge} c)\\
\hspace*{20mm}=(\beta\gamma{+}\beta\eta{+}\eta\beta{+}\alpha\alpha{+}\alpha\eta{+}\eta\alpha{+}\eta\eta)$$\\
\hspace*{30mm}{\scriptstyle\wedge}
(c'_5c''_5{+}c'_6e''_6{+}e'_7c''_7{+}a'_8a''_8{+}a'_9e''_9{+}e'_{10}a''_{10}{+}e'_{11}e''_{11})\\
\hspace*{20mm}=(\beta{\scriptstyle\wedge} c'_5)(\gamma{\scriptstyle\wedge} c''_5){+}
(\beta{\scriptstyle\wedge} c'_6)(\eta{\scriptstyle\wedge} e''_6){+}
(-1)^{|e'_7||\beta|}(\eta{\scriptstyle\wedge} e'_7)(\beta{\scriptstyle\wedge} c''_7)\\
\hspace*{23mm}{+}(\alpha{\scriptstyle\wedge} a'_8)(\alpha{\scriptstyle\wedge} a''_8){+}
(\alpha{\scriptstyle\wedge} a'_9)(\eta{\scriptstyle\wedge} e''_9){+}
(\eta{\scriptstyle\wedge} e'_{10})(\alpha{\scriptstyle\wedge} a''_{10}){+}
(\eta{\scriptstyle\wedge} e'_{11})(\eta{\scriptstyle\wedge} e''_{11})\\
\hspace*{20mm}=\phi(b'_5){\otimes}\phi(c''_5){+}\phi(b'_6){\otimes}\phi(e''_6){+}
(-1)^{|e'_7|}\phi(e'_7){\otimes}\phi(b''_7)\\
\hspace*{23mm}{+}\phi(a'_8){\otimes}\phi(a''_8){+}
\phi(a'_9){\otimes}\phi(e''_9){+}\phi(e'_{10}){\otimes}\phi(a''_{10}){+}
\phi(e'_{11}){\otimes}\phi(e''_{11})\\
\hspace*{20mm}=(\phi{\otimes}\phi)(\psi_{(X|A)}(b))$,

$\hspace*{20mm}(\psi{\otimes}_\XX\psi_{(X,A)})(\phi(a))=(\psi{\otimes}_\XX\psi_{(X,A)})(\alpha{\scriptstyle\wedge} a)\\
\hspace*{20mm}=(\alpha\alpha{+}\alpha\eta{+}\eta\alpha{+}\eta\eta){\scriptstyle\wedge}
(a'_{12}a''_{12}{+}a'_{13}e''_{13}{+}e'_{14}a''_{14}{+}e'_{15}e''_{15})\\
\hspace*{20mm}=(\alpha{\scriptstyle\wedge} a'_{12})(\alpha{\scriptstyle\wedge} a''_{12})
{+}(\alpha{\scriptstyle\wedge} a'_{13})(\eta{\scriptstyle\wedge} e''_{13})
{+}(\eta{\scriptstyle\wedge} e'_{14})(\alpha{\scriptstyle\wedge} a''_{14})
{+}(\eta{\scriptstyle\wedge} e'_{15})(\eta{\scriptstyle\wedge} e''_{15})\\
\hspace*{20mm}=\phi(a'_{12}){\otimes}\phi(a''_{12}){+}\phi(a'_{13}){\otimes}\phi(e''_{13})
{+}\phi(e'_{14}){\otimes}\phi(a''_{14}){+}\phi(e'_{15}){\otimes}\phi(e''_{15})\\
\hspace*{20mm}=(\phi{\otimes}\phi)(\psi_{(X|A)}(a))$.

Thus, $(\psi{\otimes}_\XX\psi_{(X,A)})\phi=(\phi{\otimes}\phi)\psi_{(X|A)}$.

The normal and special case is similar and easier.
\hfill$\Box$\vspace{3mm}

{\bf Definition 6.5} Let $(\underline{X},\underline{A})=\{(X_k,A_k)\}_{k=1}^m$ be a sequence of CW-pairs
such that every $(X_k,A_k)$ is homology split.

The {\it character coalgebra and algebra} of $(\underline{X},\underline{A})$ are\vspace{-2mm}
$$(C_*^{\XX_m}(\underline{X}|\underline{A}),\psi\!_{_{(\underline{X}|\underline{A})}})
=(C_*^\XX\!(X_1|A_1){\otimes}{\cdots}{\otimes}C_*^\XX\!(X_m|A_m),
\psi\!_{_{(X_1|A_1)}}{\otimes}{\cdots}{\otimes}\psi\!_{_{(X_m|A_m)}}),
\vspace{-2mm}$$
$$(C^*_{\XX_m}(\underline{X}|\underline{A}),\pi\!_{_{(\underline{X}|\underline{A})}})
=(C^*_\XX(X_1|A_1){\otimes}{\cdots}{\otimes}C^*_\XX(X_m|A_m),
\pi\!_{_{(X_1|A_1)}}{\otimes}{\cdots}{\otimes}\pi\!_{_{(X_m|A_m)}}).$$

The {\it homology coalgebra and cohomology algebra} of $(\underline{X},\underline{A})$ are\vspace{-2mm}
$$(H_*^{\XX_m}(\underline{X},\underline{A}){,}\psi\!_{_{(\underline{X},\underline{A})}})
{=}(H_*^\XX\!(X_1,A_1){\otimes}{\cdots}{\otimes}H_*^\XX\!(X_m,A_m){,}
\psi\!_{_{(X_1,A_1)}}{\otimes}{\cdots}{\otimes}\psi\!_{_{(X_m,A_m)}}),
\vspace{-2mm}$$
$$(H^*_{\XX_m}(\underline{X},\underline{A}),\pi\!_{_{(\underline{X},\underline{A})}})
{=}(H^*_\XX(X_1,A_1){\otimes}{\cdots}{\otimes}H^*_\XX(X_m,A_m){,}
\pi\!_{_{(X_1,A_1)}}{\otimes}{\cdots}{\otimes}\pi\!_{_{(X_m,A_m)}}).\vspace{3mm}$$

{\bf Theorem 6.6} {\it Let $K$ be a simplicial complex on $[m]$.

For any coproduct $\varphi\colon T_*^\XX\to T_*^\XX{\otimes}T_*^\XX$ on $T_*^\XX$
with $\varphi_m=\varphi{\otimes}{\cdots}{\otimes}\varphi$ ($m$ fold),
the total chain group $T_*^{\XX_m}(K)$ in Definition 4.5 is a subcoalgebra of $(T_*^{\XX_m},\varphi_m)$.
Denote the subcoalgebra by $(T_*^{\XX_m}(K),\varphi_K)$.

Dually, for any product $\varpi\colon T^*_\XX{\otimes}T^*_\XX\to T^*_\XX$ on $T^*_\XX$
with $\varpi_m=\varpi{\otimes}{\cdots}{\otimes}\varpi$,
the total cochain group $T^*_{\XX_m}(K)$ is a quotient algebra of $(T^*_{\XX_m},\varpi_m)$.
Denote the quotient algebra by $(T^*_{\XX_m}(K),\varpi_K)$.
\vspace{2mm}

Proof}\, $T(\tau)$ in Definition 4.5 is a subcoalgebra of $(T_*^{\XX_m},\varphi_m)$.
So $T_*^{\XX_m}(K)=+_{\tau\in K}\,T(\tau)$ is a subcoalgebra of $(T_*^{\XX_m},\varphi_m)$.
\hfill$\Box$\vspace{3mm}

{\bf Theorem 6.7} {\it The group isomorphisms in Theorem 4.6 are (co)algebra isomorphisms\vspace{-2mm}
$$(C_*^{\XX_m}(K;\underline{X},\underline{A}),\psi_{(K;\underline{X},\underline{A})})
\cong(T_*^{\XX_m}(K){\otimes}_{\XX_m}H_*^{\XX_m}(\underline{X},\underline{A}),
\psi_K{\otimes}_{\XX_m}\psi_{(\underline{X},\underline{A})}),\vspace{-2mm}$$
$$(T^*_{\XX_m}(K){\otimes}_{\XX_m}H^*_{\XX_m}(\underline{X},\underline{A}),
\pi_K{\otimes}_{\XX_m}\pi_{(\underline{X},\underline{A})})
\cong(C^*_{\XX_m}(K;\underline{X},\underline{A}),\pi_{(K;\underline{X},\underline{A})}).
\vspace{2mm}$$
If every pair $(X_k,A_k)$ is normal (or special), then $\psi_K$ and $\pi_K$
can be replaced by $\w\psi_K$ and $\w\pi_K$ (or $\overline\psi_K$ and $\overline\pi_K$) respectively.
\vspace{2mm}

\it Proof}\, By Theorem 6.4,\\
\hspace*{15mm}$(C_*^\XX\!(X_k|A_k),\psi_{(X_k|A_k)})
\cong(T_*^\XX{\otimes}_\XX H_*^\XX\!(X_k,A_k),\psi{\otimes}_\XX\psi_{(X_k,A_k)})$.\\
The restriction of the isomorphism is \\
\hspace*{22mm}$(C_*(A_k),\psi_{A_k})
\cong(S_*^\XX{\otimes}_\XX H_*^\XX\!(X_k,A_k),\psi{\otimes}_\XX\psi_{(X_k,A_k)})$.\\
By Theorem 5.6, the group isomorphism $\phi_\tau\colon H_*(\tau)\to T_*(\tau){\otimes}_{\XX_m}H_*(\underline{X},\underline{A})$
is a coalgebra isomorphism.
So $\phi_{(K;\underline{X},\underline{A})}=+_{\tau\in K}\,\phi_\tau$ is a coalgebra isomorphism.

The normal and special cases are similar.
\hfill $\Box$\vspace{3mm}

{\bf Definition 6.8} Let $K$ be a simplicial complex on $[m]$.

The product $\varpi_K$ in Theorem 6.6 induces cup product\vspace{-2mm}
$$\amalg_K\colon H^*_{\XX_m}(K)\otimes H^*_{\XX_m}(K)\to H^*_{\XX_m}(K)\vspace{-2mm}$$
defined by $[a]\amalg_K[b]=[\varpi_K(a{\otimes}b)]$ for cohomology classes $[a],[b]\in H^*_{\XX_m}(K)$.
For $\varpi_K=\pi_K,\w\pi_K,\overline\pi_K$,
$\amalg_K$ are respectively denoted by $\cup_K,\w\cup_K,\overline\cup_K$.
$H^*_{\XX_m}(K)$ with cup product $\cup_K$, $\w\cup_K$, $\overline\cup_K$ are respectively called
the {\it universal, normal, special cohomology algebra of $K$}.
\vspace{3mm}

{\bf Theorem 6.9} {\it \it For a  homology split ${\cal Z}(K;\underline{X},\underline{A})$,\vspace{-1mm}
$$\big(H^*({\cal Z}(K;\underline{X},\underline{A})),\cup\big)
\cong\big(H^*_{\XX_m}(K)\otimes_{\XX_m} H^*_{\XX_m}(\underline{X},\underline{A}),\,
\cup_{K}\,{\otimes}_{\XX_m}\,\pi_{(\underline{X},\underline{A})}\big).\vspace{-1mm}$$
If every pair $(X_k,A_k)$ is normal (or special), then
$\cup_K$ can be replaced by $\w\cup_K$ (or $\overline\cup_K$).

The conclusion holds for all polyhedral product spaces if the cohomology is taken over a field.}
\vspace{3mm}

{\it Proof}\, Corollary of Theorem 2.11 and Theorem 6.7.
\hfill$\Box$\vspace{3mm}

{\bf Example 6.10} For a homology split pair $(X,A)$, the pair $(SX,SA)$ is special, where $S$ means suspension.
So for $(\underline{SX},\underline{SA})
=\{(SX_k,SA_k)\}_{k=1}^m$ such that every $(X_k,A_k)$ is homology split,\vspace{-2mm}
$$\big(H^*({\cal Z}(K;\underline{SX},\underline{SA})),\cup\big)
\cong \big(H^*_{\XX_m}(K)\otimes_{\XX_m} H^*_{\XX_m}(\underline{SX},\underline{SA}),\,
\overline\cup_{K}\,{\otimes}_{\XX_m}\,\pi_{(\underline{SX},\underline{SA})}\big).\vspace{-2mm}$$

With the identification in Example 4.9 and regardless of the difference of even degree,
$(H^*({\cal Z}(K;S^4,S^2)),\cup)\cong (H^*_{\XX_m}(K),\overline\cup_K)$. This shows that the special cohomology algebra
of $K$ is an associative, commutative algebra.

\section{Restriction Product}

\hspace*{5mm}In this section, we will determine the restriction products of all (co)algebras defined
in Section 6.
\vspace{3mm}

{\bf Definition 7.1} Let $K$ be a simplicial complex on $[m]$.

For $(\sigma\!,\,\omega),(\sigma'\!,\,\omega'),(\sigma''\!,\,\omega'')\in\XX_m$ such that
$(\sigma'{\cup}\sigma''){\setminus}\sigma\subset\omega{\setminus}(\omega'{\cup}\omega'')$,
the {\it diagonal chain coproduct}\vspace{1mm}\\
\hspace*{20mm}$\psi_\vartriangle\,\colon
(\Sigma\w C_*(K_{\sigma,\,\omega}),d)\to(\Sigma\w C_*(K_{\sigma'\!,\,\omega'})\otimes \Sigma\w C_*(K_{\sigma''\!,\,\omega''}),d)$\vspace{1mm}\\
of $K$ is defined as follows. For $\lambda\in\Sigma\w C_*(K_{\sigma,\,\omega})$,\vspace{1mm}\\
\hspace*{49.5mm}$\psi_\vartriangle(\lambda)=\langle\mu,\nu\rangle\,
\mu\!\otimes\!\nu$,\vspace{1mm}\\
where $\mu=\lambda\cap(\omega'{\setminus}(\sigma'{\cup}\sigma''))$, $\nu=\lambda\cap((\omega''{\setminus}\omega'){\setminus}(\sigma'{\cup}\sigma''))
$\vspace{2mm},
$\langle\mu,\nu\rangle$ is the sign of the permutation $\Big(
\begin{array}{cccccc}
\scriptstyle{j_1} \!&\!\scriptstyle\!{\cdots}\!\!&\!\scriptstyle{j_u}\!&\!\scriptstyle{k_{1}} \!&\!\scriptstyle\!{\cdots}\!\!&\!\scriptstyle{k_v}\\
\scriptstyle{l_1} \!&\!\scriptstyle\!{\cdots}\!\!&\!\scriptstyle{l_u}\!&\!\scriptstyle{l_{u+1}} \!&\!\scriptstyle\!{\cdots}\!\!&\!\scriptstyle{l_{u+v}}\end{array}\Big)$
if $\mu=\{j_1,{\cdots},j_u\}$, $\nu=\{k_1,{\cdots},k_v\}$\vspace{2mm} and $\mu{\cup}\nu=\{l_1,{\cdots},l_{u+v}\}$
(all are ordered sets).

Dually, the {\it diagonal cochain product}\vspace{1mm}\\
\hspace*{20mm}$\pi_\vartriangle\,\colon
(\Sigma\w C_*(K_{\sigma'\!,\,\omega'})\otimes\Sigma\w C_*(K_{\sigma''\!,\,\omega''}),\delta)\to
(\Sigma\w C_*(K_{\sigma,\,\omega}),\delta)$\vspace{1mm}\\
of $K$ is the dual of $\psi_\vartriangle$.
Precisely, for $\mu\in\Sigma\w C^*(K_{\sigma'\!,\,\omega'})$,
$\nu\in\Sigma\w C^*(K_{\sigma''\!,\,\omega''})$,\vspace{-2mm}
$$\pi_\vartriangle(\mu{\otimes}\nu)=
\langle\mu,\nu\rangle(\Sigma_{\mu=\lambda\cap(\omega'{\setminus}(\sigma'{\cup}\sigma'')), \,\nu=\lambda\cap((\omega''{\setminus}\omega'){\setminus}(\sigma'{\cup}\sigma'')),\,
\lambda\in K_{\sigma,\omega}}\,\lambda),\vspace{-2mm}$$
where the right side is $0$ if $\mu\cap\nu\neq\emptyset$ or there is no $\lambda$ satisfying the condition.

The {\it diagonal cup product}\vspace{1mm}\\
\hspace*{24mm}$\cup_\vartriangle\,\colon
\w H^{*-1}(K_{\sigma'\!,\,\omega'})\otimes\w H^{*-1}(K_{\sigma''\!,\,\omega''})\to
\w H^{*-1}(K_{\sigma,\,\omega})$\vspace{1mm}\\
of $K$ is induced by $\pi_\vartriangle$, i.e., $[a]\cup_\vartriangle[b]=[\pi_\vartriangle(a{\otimes}b)]$ for all $[a],[b]$.
\vspace{3mm}

{\bf Theorem 7.2} {\it Let $K$ be a simplicial complex on $[m]$.
Identify $T_*^{\sigma\!,\,\omega}(K)$ and $T^*_{\sigma\!,\,\omega}(K)$
respectively with $\Sigma\w C_*(K_{\sigma\!,\,\omega})$ and $\Sigma\w C^*(K_{\sigma\!,\,\omega})$.

The restriction coproduct\vspace{-2mm}
$$\psi_{\sigma'\!,\,\omega';\sigma''\!,\,\omega''}^{\sigma\!,\,\omega},\,
\w\psi_{\sigma'\!,\,\omega';\sigma''\!,\,\omega''}^{\sigma\!,\,\omega},\,
\overline\psi_{\sigma'\!,\,\omega';\sigma''\!,\,\omega''}^{\sigma\!,\,\omega}\,\colon
T_*^{\sigma,\,\omega}(K)\to T_*^{\sigma'\!,\,\omega'}(K)\otimes T_*^{\sigma''\!,\,\omega''}(K)\vspace{-2mm}$$
of $\psi_K,\w\psi_K,\overline\psi_K$ is either the diagonal chain coproduct $\psi_\vartriangle$ or $0$,
as shown in the following table.\vspace{3mm}\\
\hspace*{3mm}{\rm \begin{tabular}{|c|c|c|}
\hline
&$=\psi_\vartriangle$&$=0$\\
\hline
{\rule[-3mm]{0mm}{8mm}$\psi_{\sigma'\!,\,\omega';\sigma''\!,\,\omega''}^{\sigma\!,\,\omega}$}
&$(\sigma'{\cup}\sigma''){\setminus}\sigma\subset\omega{\setminus}(\omega'{\cup}\omega'')$&otherwise\\
\hline
{\rule[-3mm]{0mm}{8mm}$\w\psi_{\sigma'\!,\,\omega';\sigma''\!,\,\omega''}^{\sigma\!,\,\omega}$}
& $\sigma'{\cup}\sigma''\subset\sigma$,\, $\omega\subset\omega'{\cup}\omega''$&otherwise\\
\hline
{\rule[-3mm]{0mm}{8mm}$\overline\psi_{\sigma'\!,\,\omega';\sigma''\!,\,\omega''}^{\sigma\!,\,\omega}$}
&$\sigma'{\cup}\sigma''=\sigma,\,\sigma'{\cap}\sigma''=\emptyset$,\, $\omega=\omega'{\cup}\omega'',\,\omega'{\cap}\omega''=\emptyset$&otherwise\\
\hline
\end{tabular}}\vspace{2mm}

Dually, the restriction product\vspace{-2mm}
$$\pi^{\sigma'\!,\,\omega';\sigma''\!,\,\omega''}_{\sigma\!,\,\omega},\,
\w\pi^{\sigma'\!,\,\omega';\sigma''\!,\,\omega''}_{\sigma\!,\,\omega},\,
\overline\pi^{\sigma'\!,\,\omega';\sigma''\!,\,\omega''}_{\sigma\!,\,\omega}\,\colon
T^*_{\sigma'\!,\,\omega'}(K)\otimes T^*_{\sigma''\!,\,\omega''}(K)\to T^*_{\sigma,\,\omega}(K)\vspace{-2mm}$$
of $\pi_K,\w\pi_K,\overline\pi_K$ is either the diagonal cochain product $\pi_\vartriangle$ or $0$ as follows.\vspace{3mm}\\
\hspace*{3mm}{\rm \begin{tabular}{|c|c|c|}
\hline
&$=\pi_\vartriangle$&$=0$\\
\hline
{\rule[-3mm]{0mm}{8mm}$\pi^{\sigma'\!,\,\omega';\sigma''\!,\,\omega''}_{\sigma\!,\,\omega}$}
&$(\sigma'{\cup}\sigma''){\setminus}\sigma\subset\omega{\setminus}(\omega'{\cup}\omega'')$&otherwise\\
\hline
{\rule[-3mm]{0mm}{8mm}$\w\pi^{\sigma'\!,\,\omega';\sigma''\!,\,\omega''}_{\sigma\!,\,\omega}$}
&$\sigma'{\cup}\sigma''\subset\sigma$,\, $\omega\subset\omega'{\cup}\omega''$&otherwise\\
\hline
{\rule[-3mm]{0mm}{8mm}$\overline\pi^{\sigma'\!,\,\omega';\sigma''\!,\,\omega''}_{\sigma\!,\,\omega}$}
&$\sigma'{\cup}\sigma''=\sigma,\,\sigma'{\cap}\sigma''=\emptyset$,\, $\omega=\omega'{\cup}\omega'',\,\omega'{\cap}\omega''=\emptyset$&otherwise\\
\hline
\end{tabular}}\vspace{2mm}

So the restriction product\vspace{-2mm}
$$\cup^{\sigma'\!,\,\omega';\sigma''\!,\,\omega''}_{\sigma\!,\,\omega},\,
\w\cup^{\sigma'\!,\,\omega';\sigma''\!,\,\omega''}_{\sigma\!,\,\omega},\,
\overline\cup^{\sigma'\!,\,\omega';\sigma''\!,\,\omega''}_{\sigma\!,\,\omega}\,\colon
H^*_{\sigma'\!,\,\omega'}(K)\otimes H^*_{\sigma''\!,\,\omega''}(K)\to H^*_{\sigma,\,\omega}(K)\vspace{-2mm}$$
of $\cup_K,\w\cup_K,\overline\cup_K$ is either the diagonal cup product $\cup_\vartriangle$ or $0$ as follows.\vspace{3mm}\\
\hspace*{3mm}{\rm \begin{tabular}{|c|c|c|}
\hline
&$=\cup_\vartriangle$&$=0$\\
\hline
{\rule[-3mm]{0mm}{8mm}$\cup^{\sigma'\!,\,\omega';\sigma''\!,\,\omega''}_{\sigma\!,\,\omega}$}
&$(\sigma'{\cup}\sigma''){\setminus}\sigma\subset\omega{\setminus}(\omega'{\cup}\omega'')$&otherwise\\
\hline
{\rule[-3mm]{0mm}{8mm}$\w\cup^{\sigma'\!,\,\omega';\sigma''\!,\,\omega''}_{\sigma\!,\,\omega}$}
&$\sigma'{\cup}\sigma''\subset\sigma$,\, $\omega\subset\omega'{\cup}\omega''$&otherwise\\
\hline
{\rule[-3mm]{0mm}{8mm}$\overline\cup^{\sigma'\!,\,\omega';\sigma''\!,\,\omega''}_{\sigma\!,\,\omega}$}
&$\sigma'{\cup}\sigma''=\sigma,\,\sigma'{\cap}\sigma''=\emptyset$,\, $\omega=\omega'{\cup}\omega'',\,\omega'{\cap}\omega''=\emptyset$&otherwise\\
\hline
\end{tabular}}\vspace{3mm}

Proof}\, We only prove the universal case, other cases are similar and easier.

We first prove that if $\psi_{\sigma'\!,\,\omega';\sigma''\!,\,\omega''}^{\sigma\!,\,\omega}\neq 0$,
then $(\sigma'{\cup}\sigma''){\setminus}\sigma\subset\omega{\setminus}(\omega'{\cup}\omega'')$.

Denote by $t_{W,X,Y,Z}$ the generator $t_1\otimes\cdots\otimes t_m$ of $\,T_*^{\XX_m}$
such that\\
\hspace*{7mm}$W=\{k\,|\,t_k\!=\!\alpha\}$, $X=\{k\,|\,t_k\!=\!\beta\}$, $Y=\{k\,|\,t_k\!=\!\gamma\}$, $Z=\{k\,|\,t_k\!=\!\eta\}$.\\
Suppose for $t=t_{A,B,C,D}\in T_*^{\sigma,\omega}$,\\
$\psi_m(t)=\Sigma(\pm1)
t_{A',B',C',D'}{\otimes}t_{A'',B'',C'',D''}
=\Sigma(\pm1)(t'_1{\otimes}{\cdots}{\otimes}t'_m){\otimes}(t''_1{\otimes}{\cdots}{\otimes}t''_m)$,\\
where $t_{A',B',C',D'}\in T_*^{\sigma'\!,\,\omega'}$, $t_{A'',B'',C'',D''}\in T_*^{\sigma''\!,\,\omega''}$
and $\pm1=\langle B'\!,\,B''\rangle$ with $\langle\,,\rangle$ as in Definition 7.1.
Then we have

(1) $D\subset(D'{\cup}C'){\cap}(D''{\cup}C'')$, for if $t_k=\eta$, then $t'_k,t''_k=\eta$ or $\gamma$.

(2) $C\subset C'{\cup}C''$, for if $t_k=\gamma$, then at least one of $t'_k$ and $t''_k$ is $\gamma$.

(3) $B{\setminus}(B'{\cup}B'')\subset(A'{\cup}D'){\cap}(A''{\cup}D'')$,
for if $t_k=\beta$ and $t'_k,t''_k\neq\beta$, then $t'_k,t''_k=\alpha$ or $\eta$.

(4) $A\subset(A'{\cup}D'){\cap}(A''{\cup}D'')$, for if $t_k=\alpha$, then $t'_k,t''_k=\alpha$ or $\eta$.

(5) $(B'{\cup}B'')\subset B$, for if $t'_k=\beta$ or $t''_k=\beta$, then $t_k=\beta$.

(6) $(A'{\cup}A''){\setminus}A\subset B{\setminus}(B'{\cup}B'')$, for if $t'_k=\alpha$ or $t''_k=\alpha$,
then $t_k=\alpha$ or $\beta$.

(3) implies $(B{\setminus}(B'{\cup}B'')){\cap}(C'{\cup}C'')=\emptyset$.
So (2) and (3) imply\\
\hspace*{17mm}$B{\setminus}(B'{\cup}B'')=
(B{\cup}C){\setminus}(B'{\cup}B''{\cup}C'{\cup}C'')
=\omega{\setminus}(\omega'{\cup}\omega'')$.\\
Then (6) implies $(\sigma'{\cup}\sigma''){\setminus}\sigma\subset\omega{\setminus}(\omega'{\cup}\omega'')$.

Now we prove that $\psi_{\sigma'\!,\,\omega';\sigma''\!,\,\omega''}^{\sigma\!,\,\omega}
=\psi_\vartriangle$ for $(\sigma'{\cup}\sigma''){\setminus}\sigma\subset\omega{\setminus}(\omega'{\cup}\omega'')$.

For free groups $G,G',G''$, a coproduct $\phi\colon G\to G'{\otimes}G''$ is called a base inclusion
if for every generator $g\in G$, there are unique generators $g'\in G'$ and
$g''\in G''$ such that $\phi(g)=\pm\, g'{\otimes}g''$.
It is easy to check that every restriction coproduct of $\psi\colon T_*^\XX\to T_*^\XX{\otimes}T_*^\XX$
is either a base inclusion or $0$. This implies that
every restriction coproduct $\psi_R$ of $\psi_m=\psi{\otimes}{\cdots}{\otimes}\psi$
is either a base inclusion or $0$.
So as a restriction of $\psi_R$,
$\psi_{\sigma'\!,\,\omega';\sigma''\!,\,\omega''}^{\sigma\!,\,\omega}$ (related to $K$)
is either a base inclusion or $0$.

For $(\sigma\!,\omega),(\sigma'\!\!,\omega'),(\sigma''\!\!,\omega'')\in\XX_m$ such that
$(\sigma'{\cup}\sigma''){\setminus}\sigma\subset\omega{\setminus}(\omega'{\cup}\omega'')$
and the generator $t\!=\!t_{\sigma\!,\,B,\,\omega{\setminus}B,\,[m]{\setminus}(\sigma{\cup}\omega)}\in T_*^{\sigma\!,\,\omega}(K)$ ($B\in K_{\sigma,\omega}$), $\psi_m(t)$ has a summand
$\pm\,t_{A',B',C',D'}{\otimes}t_{A'',B'',C'',D''}
=\pm(t'_1{\otimes}{\cdots}{\otimes}t'_m){\otimes}(t''_1{\otimes}{\cdots}{\otimes}t''_m)$
defined as follows.

(1) For $k\in A{\cup}\sigma'{\cup}\sigma''$, $t_k=\alpha$ or $\beta$. Take
$t'_k=\alpha$, $t''_k=\eta$ if $k\in\sigma'{\setminus}\sigma''$;
$t'_k=\eta$, $t''_k=\alpha$ if $k\in\sigma''{\setminus}\sigma'$;
$t'_k=t''_k=\alpha$ if $k\in\sigma'{\cap}\sigma''$; $t'_k=t''_k=\eta$, otherwise.

(2) For $k\in B{\setminus}(\sigma'{\cup}\sigma'')$, $t_k=\beta$.
Take $t'_k=\beta$, $t''_k=\gamma$ if $k\in\omega'{\cap}\omega''$; $t'_k=\beta$, $t''_k=\eta$ if $k\in\omega'{\setminus}\omega''$;
$t'_k=\eta$, $t''_k=\beta$ if $k\in\omega''{\setminus}\omega'$; $t'_k=t''_k=\eta$, otherwise.

(3) For $k\in C$, $t_k=\gamma$. Take
$t'_k=\gamma$, $t''_k=\eta$ if $k\in\omega'{\setminus}\omega''$;
$t'_k=\eta$, $t''_k=\gamma$ if $k\in\omega''{\setminus}\omega'$;
$t'_k=t''_k=\gamma$ if $k\in\omega'{\cap}\omega''$.

(4) For $k\in D$, $t_k=\eta$. Take
$t'_k=\gamma$, $t''_k=\eta$ if $k\in\omega'{\setminus}\omega''$;
$t'_k=\eta$, $t''_k=\gamma$ if $k\in\omega''{\setminus}\omega'$;
$t'_k=t''_k=\gamma$ if $k\in\omega'{\cap}\omega''$; $t'_k=t''_k=\eta$, otherwise.

These imply $A'=\sigma'$, $B'= B{\cap}(\omega'{\setminus}(\sigma'{\cup}\sigma''))$, $C'=\omega'{\setminus}B'$, $A''=\sigma''$, $B''= B{\cap}((\omega''{\setminus}\omega'){\setminus}(\sigma'{\cup}\sigma''))$,
$C''=\omega''{\setminus}B''$.
So as a base inclusion,\vspace{1mm}\\
\hspace*{10mm}$\psi_{\sigma'\!,\,\omega';\sigma''\!,\,\omega''}^{\sigma\!,\,\omega}
(t_{\sigma\!,\,B,\,\omega{\setminus}B,\,[m]{\setminus}(\sigma{\cap}\omega)})=
\langle B'\!,\,B''\rangle\,t_{A'\!,\,B'\!,\,C'\!,\,D'}\!\otimes\!
t_{A''\!,\,B''\!,\,\,C''\!,\,D''}$\vspace{1mm}\\
is just the $\psi_\vartriangle$ in Definition 7.1.
\hfill$\Box$\vspace{3mm}

The simplest case to compute the cohomology ring is described in the following theorem.
\vspace{2mm}

{\bf Theorem 7.3}\, {\it Let ${\cal Z}(K;\underline{X},\underline{A})$ be a homology split space
such that every\\ $i_k^*\colon H^*(X_k)\to H^*(A_k)$ is an epimorphism.
Then every ${\rm ker}\,i_k^*$ is an ideal of $H^*(X_k)$.
Define ideal $I(K)$ of $H^*(X_1){\otimes}{\cdots}{\otimes}H^*(X_m)$ by
$$I(K)=\oplus_{\sigma\notin K}\,J(\sigma),\,\,
J(\sigma)= J_1{\otimes}{\cdots}{\otimes}J_m,\,\,
J_k=\left\{\begin{array}{cl}
{\rm ker}\,i_k^*&{\rm if}\,\,k\notin\sigma,\\
H^*(A_k)&{\rm if}\,\,k\in\sigma.
\end{array}\right.\vspace{-2mm}$$
Then there is a ring isomorphism\vspace{1mm}\\
\hspace*{22mm}$H^*({\cal Z}(K;\underline{X},\underline{A}))\cong
\big(H^*(X_1){\otimes}{\cdots}{\otimes}H^*(X_m)\big)/I(K)$.\vspace{1mm}\\
Specifically, $H^*({\cal Z}(K;C\!P^\infty\!\!,*))\cong \zz(K)$,
the Stanley-Reisner ring of $K$.
}\vspace{2mm}

{\it Proof}\, Since ${\rm ker}\,i_k=0$ for all $k$,
we have $H_*^{\sigma\!,\,\omega}(\underline{X},\underline{A})=0$ if $\omega\neq\emptyset$.
Since $H_*^{\sigma,\emptyset}(K)\cong\zz$ for every $\sigma\in K$, we may identify
$H_*^{\sigma,\emptyset}(K){\otimes}H_*^{\sigma,\emptyset}(\underline{X},\underline{A})$ with
$H_*^{\sigma,\emptyset}(\underline{X},\underline{A})$. With this identification,
$H_*({\cal Z}(K;\underline{X},\underline{A}))
\cong \oplus_{\sigma\in K}\,H_*^{\sigma,\emptyset}(\underline{X},\underline{A})$
is a subcoalgebra of $H_*^{\XX_m}(\underline{X},\underline{A})$.
By definition,\vspace{1mm}\\
\hspace*{11mm}$(C_*^\XX\!(X_k|A_k),\psi_{(X_k|A_k)})=(H_*^\XX\!(X_k,A_k),\psi_{(X_k,A_k)})=(H_*(X_k),\psi_{X_k})$
\vspace{1mm}\\
with $H_*^{\emptyset,\emptyset}(X_k,A_k)= H_*(A_k)$ and
$H_*^{\{1\},\emptyset}(X_k,A_k)={\rm coker}\,i_k$.
So\vspace{1.5mm}\\
\hspace*{17mm}$H_*^{\sigma\!,\emptyset}(\underline{X},\underline{A})= H_1{\otimes}{\cdots}{\otimes}H_m,\,\,
H_k=\left\{\begin{array}{cl}
{\rm coker}\,i_k&{\rm if}\,\,k\in\sigma,\\ H_*(A_k)&{\rm if}\,\,k\notin\sigma.
\end{array}\right.$\vspace{2mm}

Dually, $H^*({\cal Z}(K;\underline{X},\underline{A}))
\cong \oplus_{\sigma\in K}\,H^*_{\sigma,\emptyset}(\underline{X},\underline{A})
\cong H^*_{\XX_m}(\underline{X},\underline{A})/I(K)$.
\hfill$\Box$\vspace{3mm}

In the remaining part,
we will compute the cohomology ring of homology split ${\cal Z}(K;\underline{X},\underline{A})$ such that
every $i_k^*\colon H^*(X_k)\to H^*(A_k)$ is a monomorphism.
In this case, $H^*_\XX(X_k,A_k)={\rm im}\,i_k^*\oplus{\rm coker}\,i_k^*= H^*(A_k)$.
So there are two products on $H^*(A_k)$.
One is the cup product $\cup_{A_k}$ induced by the diagonal map of $A_k$,
the other is $\pi_{(X_k,A_k)}$ of $H^*_\XX(X_k,A_k)$.
We always denote $H_*^\XX\!(X_k,A_k)$ and $H^*_\XX(X_k,A_k)$ respectively by $H_*(A_k)$ and $H^*(A_k)$. Then\vspace{1mm}\\
\hspace*{16.5mm}$H_*^{\XX_m}(\underline{X},\underline{A})
= H_*(A_1){\otimes}{\cdots}{\otimes}H_*(A_m)
=\oplus_{\omega\subset[m]}\,H_*^{\emptyset,\omega}(\underline{X},\underline{A})$,\vspace{1mm}\\
\hspace*{16mm}$H^*_{\XX_m}(\underline{X},\underline{A})
= H^*(A_1){\otimes}{\cdots}{\otimes}H^*(A_m)
=\oplus_{\omega\subset[m]}\,H^*_{\emptyset,\omega}(\underline{X},\underline{A})$.\vspace{1mm}\\
Regard them as groups indexed by the set\vspace{1mm}\\
\hspace*{30mm}$\RR_m=\{\,\omega\subset[m]\,\}=\{(\emptyset,\omega)\in\XX_m\}\subset\XX_m$\vspace{1mm}\\
and denote them by $H_*^{\RR_m}(\underline{X},\underline{A})$ and $H^*_{\RR_m}(\underline{X},\underline{A})$.
Then $C_*^{\XX_m}(\underline{X}|\underline{A})$ and $C^*_{\XX_m}(\underline{X}|\underline{A})$
may also be regarded as groups indexed by $\RR_m$ and denoted by
$C_*^{\RR_m}(\underline{X}|\underline{A})$ and $C^*_{\RR_m}(\underline{X}|\underline{A})$.
So \vspace{1mm}\\
\hspace*{16.5mm}$C_*^{\RR_m}(K;\underline{X},\underline{A})
\cong H_*^{\RR_m}(K){\otimes}_{\RR_m}H_*^{\RR_m}(\underline{X},\underline{A})$\\
\hspace*{43mm}$= H_*^{\RR_m}(K){\otimes}_{\RR_m}(H_*(A_1){\otimes}{\cdots}{\otimes}H_*(A_m))$,\vspace{1mm}\\
\hspace*{16mm}$C^*_{\RR_m}(K;\underline{X},\underline{A})
\cong H^*_{\RR_m}(K){\otimes}_{\RR_m}H^*_{\RR_m}(\underline{X},\underline{A})$\\
\hspace*{43mm}$= H^*_{\RR_m}(K){\otimes}_{\RR_m}(H^*(A_1){\otimes}{\cdots}{\otimes}H^*(A_m))$.
\vspace{3mm}

{\bf Definition 7.4} Define $(T_*^\RR=\zz(\beta,\gamma,\eta),d)$ to be the chain subcomplex of $(T_*^\XX,d)$
regarded as a chain complex indexed by $\RR=\{(\emptyset,\emptyset),(\emptyset,\{1\})\}\subset\XX$.

The {\it right universal coproduct} $\psi\colon(T_*^\RR,d)\to(T_*^\RR{\otimes}T_*^\RR,d)$ is defined as follows.

$\psi(\eta)=\eta{\otimes}\eta+\eta{\otimes}\gamma+\gamma{\otimes}\eta+\gamma{\otimes}\gamma$.

$\psi(\gamma)=\gamma{\otimes}\gamma+\gamma{\otimes}\eta+\eta{\otimes}\gamma$.

$\psi(\beta)=\beta{\otimes}\gamma+\beta{\otimes}\eta+\eta{\otimes}\beta+\eta{\otimes}\eta$.

The {\it right normal coproduct} $\w \psi$ is defined as follows.

$\w\psi(\eta)=\eta{\otimes}\eta+\eta{\otimes}\gamma+\gamma{\otimes}\eta+\gamma{\otimes}\gamma$.

$\w\psi(\gamma)=\gamma{\otimes}\gamma+\gamma{\otimes}\eta+\eta{\otimes}\gamma$.

$\w\psi(\beta)=\beta{\otimes}\gamma+\beta{\otimes}\eta+\eta{\otimes}\beta$.

The {\it right special coproduct} $\overline\psi$ is defined as follows.

$\overline\psi(\eta)=\eta{\otimes}\eta$.

$\overline\psi(\gamma)=\gamma{\otimes}\eta+\eta{\otimes}\gamma$.

$\overline\psi(\beta)=\beta{\otimes}\eta+\eta{\otimes}\beta$.

The {\it right strictly normal coproduct} $\hat\psi$ is defined as follows.

$\hat\psi(\eta)=\eta{\otimes}\eta$.

$\hat\psi(\gamma)=\gamma{\otimes}\gamma+\gamma{\otimes}\eta+\eta{\otimes}\gamma$.

$\hat\psi(\beta)=\beta{\otimes}\gamma+\beta{\otimes}\eta+\eta{\otimes}\beta$.

The {\it right weakly special coproduct} $\bar\psi$ is defined as follows.

$\bar\psi(\eta)=\eta{\otimes}\eta$.

$\bar\psi(\gamma)=\gamma{\otimes}\eta+\eta{\otimes}\gamma$.

$\bar\psi(\beta)=\beta{\otimes}\eta+\eta{\otimes}\beta+\eta{\otimes}\eta$.

The corresponding right coalgebras and their dual algebras are denoted as follows, \vspace{1mm}\\
\hspace*{26mm}$(T_*^{\RR_m},\varphi_m)=(T_*^\RR{\otimes}{\cdots}{\otimes}T_*^\RR,\,\varphi{\otimes}{\cdots}{\otimes}\varphi)
\,\,(m\,{\rm fold})$,\\
\hspace*{26mm}$(T^*_{\RR_m},\varpi_m)=(T^*_\RR{\otimes}{\cdots}{\otimes}T^*_\RR\,,\,\varpi{\otimes}{\cdots}{\otimes}\varpi)
\,\,(m\,{\rm fold})$,\\
where
$\varphi=\psi,\w\psi,\overline\psi,\hat\psi,\bar\psi$,
$\varpi=\pi,\w\pi,\overline\pi,\hat\pi,\bar\pi$.\vspace{3mm}

{\bf Definition 7.5} Let $K$ be a simplicial complex on $[m]$.

The {\it right total chain complex} $(T_*^{\RR_m}(K),d)$ of $K$ is
the chain subcomplex of $(T_*^{\RR_m},d)$ defined as follows.
For a subset $\tau$ of $[m]$, let\vspace{-2mm}
$$(T_*(\tau),d)=(T_1{\otimes}{\cdots}{\otimes}T_m,d),\quad
T_k=\left\{\begin{array}{cl}
T_\RR&{\rm if}\,\,k\in\tau,\vspace{1mm}\\
\zz(\gamma,\eta)&{\rm if}\,\,k\notin\tau.
\end{array}\right.\vspace{-2mm}$$
Then $(T_*^{\RR_m}(K),d)=(+_{\tau\in K}\,T_*(\tau),d)$.

Dually, the {\it right total cochain complex} $(T^*_{\RR_m}(K),\delta)$ of $K$ is
the dual of $(T_*^{\RR_m}(K),d)$.
\vspace{3mm}

{\bf Theorem 7.6} {\it Let $K$ be a simplicial complex on $[m]$.

For any coproduct $\varphi\colon T_*^\RR\to T_*^\RR{\otimes}T_*^\RR$
with $\varphi_m=\varphi{\otimes}{\cdots}{\otimes}\varphi$ ($m$ fold),
the right total chain group $T_*^{\RR_m}(K)$ is a subcoalgebra of $(T_*^{\RR_m},\varphi_m)$.
Denote the subcoalgebra by $(T_*^{\RR_m}(K),\varphi_K)$.

Dually, for any product $\varpi\colon T^*_\RR{\otimes}T^*_\RR\to T^*_\RR$
with $\varpi_m=\varpi{\otimes}{\cdots}{\otimes}\varpi$,
the right total cochain group $T^*_{\RR_m}(K)$ is a quotient algebra of $(T^*_{\RR_m},\varpi_m)$.
Denote the quotient algebra by $(T^*_{\RR_m}(K),\varpi_K)$.

Proof}\, Analogue of Theorem 6.6.\hfill$\Box$
\vspace{3mm}

{\bf Definition 7.7}
The cup product $\cup_K$, $\w\cup_K$, $\overline\cup_K$, $\hat\cup_K$, $\bar\cup_K$
of $H^*_{\RR_m}(K)=H^*(T^*_{\RR_m}(K))$
induced by  $\pi_K$, $\w\pi_K$, $\overline\pi_K$, $\hat\pi_K$, $\bar\pi_K$
are respectively called {\it the right universal, normal, special, strictly normal, weakly special
cohomology algebra of $K$.} \vspace{3mm}

{\bf Theorem 7.8} {\it The restriction products of the right cohomology algebras of $K$ satisfy the following table.
\vspace{3mm}\\
\hspace*{22mm}
{\rm \begin{tabular}{|c|c|c|}
\hline
&$=\cup_\vartriangle$&$=0$\\
\hline
{\rule[-3mm]{0mm}{8mm}$(\cup_K)^{\omega';\omega''}_{\omega}$}
&all&\\
\hline
{\rule[-3mm]{0mm}{8mm}$(\w\cup_K)^{\omega';\omega''}_{\omega}$}
&$\omega\subset\omega'{\cup}\,\omega''$&otherwise\\
\hline
{\rule[-3mm]{0mm}{8mm}$(\overline\cup_K)^{\omega';\omega''}_{\omega}$}
& $\omega=\omega'{\cup}\,\omega''$,\,\,$\omega'{\cap}\omega''=\emptyset$&otherwise\\
\hline
{\rule[-3mm]{0mm}{8mm}$(\hat\cup_K)^{\omega';\omega''}_{\omega}$}
&$\omega=\omega'{\cup}\,\omega''$&otherwise\\
\hline
{\rule[-3mm]{0mm}{8mm}$(\bar\cup_K)^{\omega';\omega''}_{\omega}$}
& $\omega'{\cup}\,\omega''\subset\omega$,\,\,$\omega'{\cap}\omega''=\emptyset$&otherwise\\
\hline
\end{tabular}}\vspace{2mm}

Proof}\, Analogue of Theorem 7.2.
\hfill$\Box$
\vspace{6mm}

{\bf Example 7.9} Let $K$ be an non-empty complex with $m$-vertices. Then by Hochster theorem,\vspace{-2mm}
$$H^*_{\RR_m}(K)=\oplus_{\omega\in\RR_m}\w H^{*-1}(K|_\omega)
= {\rm Tor}_{\zz[x_1,{\cdots},x_m]}^*(\zz(K),\zz),\vspace{-2mm}$$
where $\zz(K)$ is the Stanley-Reisner ring of $K$.

If $K$ is the $m$-gon, $m>2$, then the vertex set is $[m]$ with edges $\{i,i{+}1\}$ for $i\in\zz_{m}$,
where $\zz_m$ is the group of integers modular $m$ regarded only as a set.
The non-zero cohomology groups $H^*_{\omega}=\w H^{*-1}(K|_\omega)$
are as follows.

(1) $H^0_{\emptyset}=\zz$ with the generator $1$ represented by
$\emptyset\in\Sigma\w C^{-1}(K_{\emptyset,\emptyset})$.

(2) Let $\omega$ be a subset of $[m]$ with connected component $\omega_1,\cdots,\omega_k$ ($k>1$).
Then $\Sigma_{u\in\omega_i}\{u\}\in\Sigma\w C^0(K_{\emptyset,\omega})$ represents a cohomology class in $H^1_{\omega}$
which is denoted by $[\omega|\omega_i]$.
$H^1_{\omega}$ is the group generated by $[\omega|\omega_1],\cdots,[\omega|\omega_k]$ modulo the zero relation $\Sigma_{i=1}^k[\omega|\omega_i]=0$.

(3) $H^2_{[m]}=\zz$ with the generator $\kappa$ represented by any directed edge $\{i,i{+}1\}
\in\Sigma\w C^1(K)$.

Denote the diagonal cup product $\cup_\vartriangle$ by ${\cup}_{\,\omega}^{\,\omega'\!,\,\omega''}\colon
H^*_{\omega'}{\otimes}H^*_{\omega''}\to H^*_\omega$.
Since $H^2_{\omega}=0$ except $\omega=[m]$,
we have $[\omega'|\omega'_i]{\cup}_{\,\omega}^{\,\omega'\!,\,\omega''}[\omega''|\omega''_j]=0$ if $\omega\neq[m]$.
So $\cup_K=\Sigma_\omega\,{\cup}_{\,\omega}^{\,\omega'\!,\,\omega''}={\cup}_{[m]}^{\,\omega'\!,\,\omega''}$.
Define product $*\colon\zz_m\times\zz_m\to\zz$ by\vspace{-2mm}
$$i*j=\left\{\begin{array}{rl}
1&{\rm if}\,\, j\equiv i{+}1\,\,{\rm mod}\,m,\\
-1&{\rm if}\,\, j\equiv i{-}1\,\,{\rm mod}\,m,\\
0&{\rm otherwise}.
             \end{array}
\right.\vspace{-2mm}$$
For subsets $A,B$ of $\zz_m$, define $A*B=\Sigma_{i\in A,\,j\in B}\,i*j$.
Then,\vspace{1mm}\\
\hspace*{20mm}$[\omega'|\omega'_i]\,{\cup}_K\,[\omega''|\omega''_j]=
[\omega'|\omega'_i]\,{\cup}_{\,[m]}^{\,\omega'\!,\,\omega''}[\omega''|\omega''_j]=
(\omega'_i*\omega''_j)\kappa$.
\vspace{3mm}

{\bf Definition 7.10}\, Let\, $(X,A)$\, be a homology split pair such that every
$i_*\colon H_*(A)\to H_*(X)$ is an epimorphism.

$(X,A)$ is called {\it right normal} if the character coproduct satisfies

$\psi_{(X|A)}(\eta)\subset\eta{\otimes}\eta+\eta{\otimes}\gamma+\gamma{\otimes}\eta+\gamma{\otimes}\gamma$.

$\psi_{(X|A)}(\gamma)\subset\gamma{\otimes}\gamma+\gamma{\otimes}\eta+\eta{\otimes}\gamma$.

$\psi_{(X|A)}(\beta)\subset\beta{\otimes}\gamma+\beta{\otimes}\eta+\eta{\otimes}\beta$.

$(X,A)$ is called {\it right special} if the character coproduct satisfies

$\psi_{(X|A)}(\eta)\subset\eta{\otimes}\eta$.

$\psi_{(X|A)}(\gamma)\subset\gamma{\otimes}\eta+\eta{\otimes}\gamma$.

$\psi_{(X|A)}(\beta)\subset\beta{\otimes}\eta+\eta{\otimes}\beta$.

$(X,A)$ is called {\it right strictly normal} if the character coproduct satisfies

$\psi_{(X|A)}(\eta)\subset\eta{\otimes}\eta$.

$\psi_{(X|A)}(\gamma)\subset\gamma{\otimes}\gamma+\gamma{\otimes}\eta+\eta{\otimes}\gamma$.

$\psi_{(X|A)}(\beta)\subset\beta{\otimes}\gamma+\beta{\otimes}\eta+\eta{\otimes}\beta$.

$(X,A)$ is called {\it right weakly special} if the character coproduct satisfies

$\psi_{(X|A)}(\eta)\subset\eta{\otimes}\eta$.

$\psi_{(X|A)}(\gamma)\subset\gamma{\otimes}\eta+\eta{\otimes}\gamma$.

$\psi_{(X|A)}(\beta)\subset\beta{\otimes}\eta+\eta{\otimes}\beta+\eta{\otimes}\eta$.
\vspace{3mm}

{\bf Theorem 7.11}\, {\it Let ${\cal Z}(K;\underline{X},\underline{A})$ be a homology split space
such that every $i_k^*\colon H^*(X_k)\to H^*(A_k)$ is a monomorphism.
Then\vspace{-2mm}
$$\big(H^*({\cal Z}(K;\underline{X},\underline{A})),\cup\big)\hspace{90mm}\vspace{-2mm}$$
$$\cong\Big(H^*_{\RR_m}(K){\otimes}_{\RR_m}(H^*(A_1){\otimes}{\cdots}{\otimes}H^*(A_m)),\,
\cup_K{\otimes}_{\RR_m}(\pi\!_{_{(X_1,A_1)}}\otimes\cdots\otimes\pi\!_{_{(X_m,A_m)}})\Big).$$

If every pair $(X_k,A_k)$ is right normal (or strictly normal, special, weakly special),
then $\cup_K$ can be replaced by $\w\cup_K$
(or $\overline\cup_K$, $\hat\cup_K$, $\bar\cup_K$).

Proof}\, Analogue of Theorem 6.9.\hfill$\Box$
\vspace{3mm}

{\bf Example 7.12} Let ${\cal Z}(K;\underline{CX},\underline{X})$, $(\underline{CX},\underline{X})=\{(CX_k,X_k)\}_{k=1}^m$,
be a polyhedral product space such that every $H_*(X_k)$ is a free group,
where $C$ means the cone of a CW-complex. Then\vspace{1mm}\\
${\rm im}\,i_k=\zz(1),\,\,
{\rm ker}\,i_k=\w H_*(X_k),\,\,
C_*^\RR(CX_k|X_k)=\w H_*(X_k)\oplus\Sigma\w H_*(X_k)\oplus\zz(1)$,\vspace{1mm}\\
where $1$ is represented by the base point.
For $a\in\w H_*(X_k)$, suppose $\psi_{X_k}(a)=1{\otimes}a+a{\otimes}1+\Sigma\, a'_i{\otimes}a''_i$
with $a'_i,a''_i\in\w H_*(X_k)$. Then
$\psi_{(CX_k|X_k)}(\overline a )=
1{\otimes}\overline a +\overline a {\otimes}1+\Sigma\,\overline a '_i{\otimes}a''_i$,
where $d\,\overline x= x$.
So $(CX_k,X_k)$ is right strictly normal and\vspace{-2mm}
$$\big(H_*({\cal Z}(K;\underline{CX},\underline{X})),\cup\big)\hspace{71mm}\vspace{-2mm}$$
$$\cong\Big(H_*^{\RR_m}(K){\otimes}_{\RR_m}(H_*(X_1){\otimes}{\cdots}{\otimes}H_*(X_m)),\,
\hat\cup_K{\otimes}_{\RR_m}(\cup_{X_1}{\otimes}{\cdots}{\otimes}\cup_{X_m})\Big),$$
where $a_1{\otimes}{\cdots}{\otimes}a_m\in H^*_{\emptyset,\omega}(\underline{CX},\underline{X})$
with $\omega=\{k\,|\,a_k\neq 1\}$.

Similarly, $(CSX_k,SX_k)$ is right special and we have\vspace{-2mm}
$$\big(H^*({\cal Z}(K;\underline{CSX},\underline{SX})),\cup\big)\hspace{71mm}\vspace{-2mm}$$
$$\cong
\Big(H^*_{\RR_m}(K){\otimes}_{\RR_m}(H^*(SX_1){\otimes}{\cdots}{\otimes}H^*(SX_m)),\,
\overline\cup_K{\otimes}_{\RR_m}(\cup_{SX_1}{\otimes}{\cdots}{\otimes}\cup_{SX_m})\Big).$$

Since the generalized moment-angle complex ${\cal Z}(K;D^n,S^{n-1})$ satisfies that
$H^*_{\emptyset,\omega}(\underline{D^n},\underline{S^{n-1}})\cong\zz$ for all $\omega$, we may identify
$H^*_{\RR_m}(K){\otimes}_{\RR_m}H^*_{\RR_m}(\underline{D^n},\underline{S^{n-1}})$ with $H^*_{\RR_m}(K)$
(degree uplifted).
With this identification,\vspace{1mm}\\
\hspace*{31mm}$\big(H^*({\cal Z}(K;D^1,S^0)),\cup\big)\cong
\big(H^*_{\RR_m}(K),\hat\cup_K\big)$,\vspace{0mm}\\
\hspace*{31mm}$\big(H^*({\cal Z}(K;D^2,S^1)),\cup\big)\cong
\big(H^*_{\RR_m}(K),\overline\cup'_K\big)$,\\
\hspace*{31mm}$\big(H^*({\cal Z}(K;D^3,S^2)),\cup\big)\cong
\big(H^*_{\RR_m}(K),\overline\cup_K\big)$,\\
where $\overline\cup'_K$ satisfies
$a\,\overline\cup'_K\,b=\langle\omega'\!,\omega''\rangle(a\,\overline\cup_K\,b)$
for $a\in H^*_{\emptyset,\omega'}(K)$,
$b\in H^*_{\emptyset,\omega''}(K)$
with $\langle\,,\rangle$ as in Definition 7.1.

These examples show that the corresponding right cohomology algebras of $K$ is an associative, commutative algebra
with unit.

With the notations in Example 7.9, we compute the cohomology ring
$H^*({\cal Z}(K;D^2,S^1))$ for $m=6$ and $K$ being the hexagon.
The generators of $H^k({\cal Z}(K;D^2,S^1))$ for $0<k<8$ are the following.

$[i,i{+}2\,|\,i]{+}[i,i{+}2\,|\,i{+}2]=0$, $i=1,{\cdots},6$,

$[j,j{+}3\,|\,j]{+}[j,j{+}3\,|\,j{+}3]=0$, $j=1,2,3$,

$[i,i{+}1,i{+}3\,|\,i,i{+}1]{+}[i,i{+}1,i{+}3\,|\,i{+}3]=0$, $i=1,{\cdots},6$,

$[i,i{+}1,i{+}4\,|\,i,i{+}1]{+}[i,i{+}1,i{+}4\,|\,i{+}4]=0$, $i=1,{\cdots},6$,

$[1,3,5\,|\,1]{+}[1,3,5\,|\,3]{+}[1,3,5\,|\,5]=0$,

$[2,4,6\,|\,2]{+}[2,4,6\,|\,4]{+}[2,4,6\,|\,6]=0$,

$[i,i{+}1,i{+}2,i{+}4\,|\,i,i{+}1,i{+}2]{+}[i,i{+}1,i{+}2,i{+}4\,|\,i{+}4]=0$, $i=1,{\cdots},6$,

$[j,j{+}1,j{+}3,j{+}4\,|\,j,j{+}1]{+}[j,j{+}1,j{+}3,j{+}4\,|\,j{+}3,j{+}4]=0$, $j=1,2,3$.

So all the non-zero cup products $H^1_{\omega'}\cup H^1_{\omega''}$ are the following.

$[i,i{+}2\,|\,i]\,\overline\cup'_K\,[i{+}3,i{+}4,i{+}5,i{+}1\,|\,i{+}1]$, $i=1,{\cdots},6$,

$[j,j{+}3\,|\,j]\,\overline\cup'_K\,[j{+}1,j{+}2,j{+}4,j{+}5\,|\,j{+}1,j{+}2]$, $j=1,2,3$,

$[i,i{+}1,i{+}3\,|\,i{+}3]\,\overline\cup'_K\,[i{+}4,i{+}5,i{+}2\,|\,i{+}2]$, $i=1,{\cdots},6$,

$[1,3,5\,|\,1]\,\overline\cup'_K\,[2,4,6\,|\,2]$, $[1,3,5\,|\,3]\,\overline\cup'_K\,[2,4,6\,|\,4]$.

This is in accordance with the fact that ${\cal Z}(K;D^2,S^1)$ is homotopic equivalent to
the connected sum of $9$ copies of $S^3{\times}S^5$ and $8$ copies of $S^4{\times}S^4$.
\vspace{3mm}

We finish the paper with an example which shows that the ring structure of
$H^*({\cal Z}(K;\underline{X},\underline{A}))$
depends not only on all $i_k^*\colon H^*(X_k)\to H^*(A_k)$,
but also the character coproducts of $(X_k,A_k)$.
\vspace{2mm}

{\bf Example 7.13}\, Let $(X_i,A_i)$ be as in Example 2.9. By definition,\vspace{1mm}\\
\hspace*{18mm}$(C_\RR(X_1|A_1),\psi_{(X_1|A_1)})\cong
(T_\RR\otimes_\RR H_*(A_1),\overline\psi\otimes_\RR\psi_{A_1})$,\\
\hspace*{18mm}$(C_\RR(X_2|A_2),\psi_{(X_2|A_2)})\cong
(T_\RR\otimes_\RR H_*(A_2),\bar\psi\otimes_\RR\psi_{(X_2,A_2)})$,\vspace{1mm}\\
where $\overline\psi$ and $\bar\psi$ are as in Definition 7.4. So\vspace{1mm}\\
\hspace*{7mm}$(H^*({\cal Z}(K;X_1,A_1)),\cup)\cong\Big(H^*_{\RR_m}(K){\otimes}_{\RR_m}H^*(A_1)^{\otimes m},
\overline\cup_K{\otimes}_{\RR_m}\cup_{A_1}^{\otimes m}\Big)$,\vspace{1mm}\\
\hspace*{7mm}$(H^*({\cal Z}(K;X_2,A_2)),\cup)\cong\Big(H^*_{\RR_m}(K){\otimes}_{\RR_m}H^*(A_2)^{\otimes m},
\bar\cup_K{\otimes}_{\RR_m}\pi_{(X_2,A_2)}^{\otimes m}\Big)$.\vspace{1mm}

$(H^*(A_1),\cup_{A_1})=\zz[x,y]/(x^2,y^2,xy)$,\,\, where $\zz[-]$ means the polynomial
algebra, the unit $1$ of $\zz[x,y]/(x^2,y^2,xy)$ is a generator of $H^0(A_1)$,
$x$ is a generator of $H^2(A_1)$ and $y$ is a generator of $H^3(A_1)$. So\vspace{1mm}\\
$(H^*(A_1)^{\otimes m},\cup_{A_1}^{\otimes m})=\zz[x_1,{\cdots},x_m,y_1,{\cdots},y_m]/(x_i^2,y_i^2,x_iy_i)
=\oplus\,H^*_{\emptyset,\omega}(\underline{X_1},\underline{A_1})$,\vspace{1mm}
where $x_{i_1}{\cdots}x_{i_s}y_{j_1}{\cdots}y_{j_t}\in
H^*_{\emptyset,\{j_1,{\cdots},j_t\}}(\underline{X_1},\underline{A_1})$.\vspace{1mm}

$(H^*(A_2),\pi_{(X_2,A_2)})=\zz[x,y]/(x^2{-}y,y^2,xy)$ (the product does not keep degree!) and so\vspace{1mm}\\
\hspace*{10mm}$(H^*(A_2)^{\otimes m},\pi_{(X_2,A_2)}^{\otimes m})
=\zz[x_1,{\cdots},x_m,y_1,{\cdots},y_m]/(x_i^2{-}y_i,y_i^2,x_iy_i)$,\vspace{1mm}\\
where $x_{i_1}{\cdots}x_{i_s}y_{j_1}{\cdots}y_{j_t}\in
H^*_{\emptyset,\{j_1,{\cdots},j_t\}}(\underline{X_2},\underline{A_2})$.
So we have ring isomorphisms
$$H^*({\cal Z}(K;X_1,A_1))\cong H^*_{\RR_m}(K){\otimes}_{\RR_m}
\zz[x_1,{\cdots},x_m,y_1,{\cdots},y_m]/(x_i^2,y_i^2,x_iy_i),\quad\,\,\vspace{-2mm}$$
$$H^*({\cal Z}(K;X_2,A_2))\cong H^*_{\RR_m}(K){\otimes}_{\RR_m}
\zz[x_1,{\cdots},x_m,y_1,{\cdots},y_m]/(x_i^2{-}y_i,y_i^2,x_iy_i),$$
where the upper $H^*_{\RR_m}(K)$ is a right special algebra and the lower $H^*_{\RR_m}(K)$
is a right weakly special algebra, although both can be right universal algebras.
\vspace{1mm}

{\bf Acknowledgment} I express my deepest gratitude to Professor Anthony Bahri for his very instructive suggestions.
I must also thank the referee for correcting so many errors and giving so helpful suggestions
to my submissions.

\end{document}